\documentclass[leqno]{article}
\usepackage[T1]{fontenc}
\usepackage[utf8]{inputenc}
\usepackage{palatino}
\usepackage[english]{babel}
\usepackage{amssymb, amsfonts,amsmath,amsthm,enumitem,array,mathtools}
\usepackage{mathrsfs}
\usepackage{tikz-cd}
\usepackage{stmaryrd}
\usepackage{stmaryrd}
\usepackage{faktor} %per gli insiemi quoziente
\usepackage{dsfont}
\usepackage{indentfirst}
\usepackage{mathdots}
\usepackage{pgfplots}
\pgfplotsset{compat=1.18}
\usepackage{tikz}
\usetikzlibrary{shapes}
\usetikzlibrary{plotmarks}
\usepackage{graphics}
\usepackage{bigints}
\usepackage{xcolor}
\usepackage{comment}    
\usepackage{csquotes}
\usepackage{hyperref}
\usepackage{mathabx}

\usepackage{titlesec}
\titleformat*{\section}{\normalfont\scshape\MakeUppercase}
\titleformat*{\subsection}{\normalfont\scshape}
\titleformat*{\subsubsection}{\normalfont\small\scshape}

\renewcommand{\thesection}{\Roman{section}}
\renewcommand{\thesubsection}{\Roman{section}.\small\arabic{subsection}}

\titleformat{\section}
  {\large\normalfont\scshape\centering}
  {\thesection.}{0.5em}{}

\titleformat{\subsection}[runin]
  {\normalfont\scshape}
  {\thesubsection.}{0.5em}{}[.]

\usepackage{tocloft}

\newtheoremstyle{customremark}% name
{10pt}% Space above
{10pt}% Space below
{\normalfont}% Body font
{}% Indent amount
{\scshape}% Theorem head font
{.}% Punctuation after theorem head
{.5em}% Space after theorem head
{}% Theorem head spec (can be left empty, meaning ‘normal’)

\newtheoremstyle{customth}
{10pt}{10pt}{\itshape}{}{\small\scshape}{.}{.5em}{}
                
\newtheoremstyle{customdef}
{10pt}{10pt}{\normalfont}{}{\small\scshape}{.}{.5em}{}

\theoremstyle{customth}
\newtheorem{thm}{Theorem} % Creates the 'thm' environment and 'thm' counter
\newtheorem{lem}[thm]{Lemma}
\newtheorem{cor}[thm]{Corollary}
\newtheorem{prop}[thm]{Proposition}

\newtheorem*{thm*}{Theorem}

\theoremstyle{customremark}
\newtheorem{oss}[thm]{Remark}

\newtheorem*{notazz}{Notations}

\theoremstyle{customdef}
\newtheorem{defn}[thm]{Definition}

\setenumerate{label*=\arabic*.}

\usepackage{geometry}
\makeatletter
\newenvironment{myproof}[1][Proof]{\par
    \pushQED{\qed}%
    \normalfont \topsep6\p@\@plus6\p@\relax
    \trivlist
    \item[\hskip\labelsep
        \small\scshape
        #1\@addpunct{.} ]\mbox{}\nobreak}
    {\popQED\endtrivlist\@endpefalse}
\makeatother
\usepackage{lipsum}

\definecolor{DarkRed}{HTML}{B6321C}

\newcommand*\xbar[1]{%
  \hbox{%
    \vbox{%
      \hrule height 0.5pt % The actual bar
      \kern0.5ex%         % Distance between bar and symbol
      \hbox{%
        \kern-0.4em%      % Shortening on the left side
        \ensuremath{#1}%
        \kern-0.4em%      % Shortening on the right side
      }%
    }%
  }%
} 

\newcommand{\xdasharrow}[2][->]{
\tikz[baseline=-\the\dimexpr\fontdimen22\textfont2\relax]{
\node[anchor=south,font=\scriptsize, inner ysep=1.5pt,outer xsep=2.2pt](x){#2};
\draw[shorten <=3.4pt,shorten >=3.4pt,dashed,#1](x.south west)--(x.south east);
}
}

\newcommand{\Ll}{\mathrm{L}}

\newcommand{\Lip}{\mathrm{Lip}}

\newcommand{\spann}{\mathrm{span\,}}

\newcommand{\sgn}{\mathrm{sgn}}

\newcommand{\Bexp}{\mathfrak{B}_{\exp}}

\newcommand{\Id}{\mathrm{Id}}

\newcommand{\C}{\mathbb{C}}
\newcommand{\R}{\mathbb{R}}

\newcommand{\N}{\mathbb{N}}

\newcommand{\cl}{\mathcal{L}}
\newcommand{\B}{\mathcal{B}}

\newcommand{\eps}{\varepsilon}
\newcommand{\de}{\partial}

\newcommand{\funcdef}[4]{
\begin{array}{ccc}
    {#1} & \to & {#2}\\
     {#3} & \mapsto & {#4}
\end{array}}

\newcommand{\funcDef}[5]{
\begin{array}{cccc}
{#1:} & {#2} & \to & {#3}\\
     & {#4} & \mapsto & {#5}
\end{array}}

\title{\normalfont{\scshape{Extreme Value Theory for discrete time samplings of  weakly dissipative SDEs by classical operator perturbation theory}}}

\author {\scshape{Antonio Marinucci}\thanks{University of Rome Tor Vergata, Italy. Email: antonio.marinucci@axp.mat.uniroma2.it} \and 
\scshape{Stefano Galatolo}\thanks{University of Pisa, Italy. Email: stefano.galatolo@unipi.it}}

\date{\scshape{\today}}

\usepackage{thmtools}
\usepackage{thm-restate}

\begin{document}

\maketitle

% --- Customize the Abstract Title ---
% \normalfont strips the hardcoded bold formatting, and we write it in all caps
\renewcommand{\abstractname}{\normalfont\scshape Abstract}

\begin{abstract}
We establish an extreme value law for discrete-time observations
of stochastic differential equations (SDEs) on $\mathbb{R}^d$ driven
by additive Brownian noise.
We consider bounded, globally Lipschitz drifts satisfying
$\langle b(x),x\rangle\leq-\kappa|x|$ outside a bounded region,
with $\kappa>0$.
For stationary initial data and any fixed sampling interval
$h>0$, we prove that the probability of avoiding a shrinking
target ball $B_n$ during the first $n$ observations converges
to $e^{-\tau}$ whenever $n\mu(B_n)\to\tau>0$, where $\mu$
is the invariant probability measure.
Equivalently, logarithmic distance observables satisfy an
extreme value law with extremal index one.
Our approach exploits the smoothing action of the noise
to apply classical spectral perturbation theory on a Banach
space of continuously differentiable densities with
exponentially weighted integrability.
On this space, the unperturbed transfer operator has a
spectral gap, and the killed transfer operators, which
remove mass in the target before evolving the remaining
density, form a continuously differentiable family in
operator norm when parameterized by target volume.
An exponential Lyapunov estimate and bounds on the transition
density and its spatial derivatives provide the required
spectral and regularity properties.
The extreme value law follows from the resulting first-order
expansion of the leading eigenvalue.
\end{abstract}
%\begin{abstract}
%    We investigate the distribution of extreme events for stochastic processes arising from the discrete-time sampling of Stochastic Differential Equations (SDEs) on $\mathbb{R}^{d}$. By employing a suitable functional analytic framework, we study the spectral properties of perturbed and unperturbed transfer operators by the classical operator perturbation theory of Hennion and Hervé.
%\end{abstract}



% --- Customize the Table of Contents Title ---
% Removes bold and sets the size (adjust \Large if you want it bigger/smaller)
\renewcommand{\cfttoctitlefont}{\normalfont\Large\scshape} 
% Changes the actual text to all caps
\renewcommand{\contentsname}{Contents}

\setcounter{tocdepth}{2}
\tableofcontents

\section{Introduction}\label{sec: intro}

Extreme value theory studies the asymptotic behaviour of large
observations along a stochastic or dynamical process. For dependent observations,
the limiting laws reflect both the frequency of rare events and
their tendency to occur in clusters
\cite{leadbetter1988extremal,leadbetterVecchio,hsing1988exceedance}.
In dynamical systems, an important connection relates extreme values
of observables that become unbounded near a prescribed point to
hitting times of shrinking neighbourhoods of that point
\cite{hirata1999statistics,haydn2005hitting,Hitting}.
This connection provides a common framework for studying extremes
and recurrence, with applications to the statistical description
of complex systems \cite{lucarini2016extremes}.
In this paper, we consider discrete-time observations of diffusion
processes and investigate how regularization by noise can be used
to analyse the associated shrinking-target problem through
classical spectral perturbation theory.

We study stochastic differential equations on $\R^d$ of the form
\[
    dX_t=b(X_t)\,dt+dW_t,
\]
where $W_t$ is a standard $d$-dimensional Brownian motion and the
drift $b$ is bounded, globally Lipschitz continuous, and satisfies
an inward-pointing condition outside a bounded region.
We consider the stationary diffusion sampled at a fixed positive
time interval and study its visits to shrinking neighbourhoods
of a prescribed point.
We prove an exponential law for the suitably rescaled first
hitting times, or equivalently an extreme value law for an
observable measuring proximity to the target on a logarithmic
scale.
The precise assumptions and statement are given in Section~2.

Exponential hitting-time laws have been established in several
general probabilistic settings.
Regenerative methods yield such limits for rare sets in positive
recurrent Harris chains and processes \cite{Glynn11}.
Mixing-based approaches provide exponential approximations and
limit laws for rare events, including events specified by long
strings of observations
\cite{abadi2001inequalities,AbadiSaussol11}.

Poisson and compound Poisson laws for visit counts to shrinking
targets have also been established for mixing dynamical systems,
with the compound Poisson case accounting for clustering near
periodic points \cite{HV04,nicolai2009compound}.
More recently, compound Poisson limits in random dynamical
systems have been obtained through spectral perturbation
methods \cite{AFGTV25} and probabilistic block approximations
\cite{AHV25}.
Abstract criteria for Poisson limits and joint limits of hitting
times and positions are provided in \cite{Zwe22}.
These works place the exponential law within a broad theory of
rare events.
In applying such results to a particular process and family of
targets, the identification of the normalization and the control
of short returns remain essential.

A particularly fruitful approach to rare-event statistics uses
transfer operators to express avoidance probabilities in terms
of the evolution of densities.
The target is regarded as a hole, and trajectories are removed
when they enter it.
The probability of avoiding the target is then represented by
iterates of the transfer operator of an open system, whose leading eigenvalue
determines the asymptotic decay.
Keller and Liverani \cite{KL09} established first-order
perturbation formulae relating the leading eigenvalue to rare
events and escape rates.
Keller \cite{KL12} subsequently developed the connection with
exponential hitting-time laws and extremal indices.
These results are part of a broader use of spectral methods
for limit theorems in dynamical systems and Markov processes
\cite{guivarc1988theoremes,HennionHerve,KM03,HP10}.

Introducing a hole involves multiplication by an indicator
function and therefore creates a difficulty for perturbation
theory on spaces of regular functions.
Even a small hole can produce a discontinuity, so smallness of
the target does not by itself imply smallness of the perturbation
in the operator norm of a given regularity space.
The Keller--Liverani theorem \cite{KL99} addresses this type of
difficulty by combining estimates in a strong norm with
perturbative control in a weaker norm.
In \cite{EXT}, the authors  used
this framework to prove an extreme value law and Poisson
statistics for discrete-time samplings of dissipative SDEs.
Their construction relies on weighted spaces of bounded
variation, which accommodate the discontinuities generated
by the holes.

The present work considers a different regime of confinement
and develops a perturbative construction on a space of regular
densities.
In \cite{EXT}, the globally Lipschitz drift satisfies a
quadratic dissipativity condition of the form
$\langle b(x),x\rangle\leq R_1-R_2|x|^2$, with $R_2>0$.
Here we assume instead that
$\langle b(x),x\rangle\leq-\kappa|x|$ outside a bounded region,
together with boundedness of $b$.
The radial confinement is therefore weaker, while boundedness
is an additional restriction.
The two results concern complementary classes of diffusions:
for example, the drift
$b(x)=-a x/\sqrt{1+|x|^2}$, with $a>0$, satisfies our assumptions
but not the quadratic dissipativity condition, whereas the
linear drift $b(x)=-x$ is covered by \cite{EXT} and is excluded
here.
Despite the weaker confinement, an exponential Lyapunov function
provides geometric control of the dynamics at infinity.
The general connection between Lyapunov conditions, local
minorization and geometric ergodicity is described, for example,
in \cite{HM11}.

The main contribution of this work is to establish a regular
functional framework in which the smoothing effect of additive
noise makes shrinking-target perturbations amenable to classical
spectral perturbation theory.
For bounded, weakly confining drifts, we work on a Banach space
of continuously differentiable densities with exponentially
weighted integrability.
On this space, we establish a spectral gap for the unperturbed
transfer operator and prove that the family of killed operators,
parameterized by target volume, is continuously differentiable
in operator norm.

This construction combines the control at infinity provided
by an exponential Lyapunov estimate with the regularizing
properties of the transition kernel.
The latter rely on the density and gradient bounds of
\cite{MPZ} and on a mixed spatial derivative estimate
established in Appendix~\ref{app: StimaD derivata mista}.

Classical spectral perturbation theory \cite{Kato}, in the
formulation of Hennion--Herv\'e \cite{HennionHerve}, then
provides a first-order expansion of the leading eigenvalue:
to leading order, its displacement from one equals the
invariant measure of the target.
Together with the spectral decomposition, this expansion
yields the extreme value law.

Throughout the paper, the sampling interval $h$ is fixed and
strictly positive.
Although the limiting law has the same form for every such $h$,
the estimates used in the proof depend on $h$.
No uniformity as $h\to0$ is asserted, and passage to continuous
time requires a separate analysis; see
\cite[Appendix~A]{EXT}.

\section{Main result}\label{sec: mainresults}

We consider a system whose continuous-time evolution is described by an SDE defined on $\R^d$ of the following type:
\begin{equation}\label{eq:system1}
    \begin{cases}
        dX_t=b(X_t)\, dt+dW_{t}\\
        X_0=x 
    \end{cases}\ ,
\end{equation}
where $W_t$ is the Brownian motion defined on a probability space $(\Omega,\mathcal{F},\mathbb{P})$. We make the following assumptions on the drift function $b:\R^d\to\R^d$:
\begin{itemize}
    \item[$(H)$] \textit{boundedness and Lipschitz continuity:} $b\in L^\infty(\R^d)$ and there exists $K>0$ such that for any $x,y\in\R^d$,
    \[
        |b(x)-b(y)|\leq K|x-y|\ .
    \]
    \item [$(AD)$] \textit{dissipativity:} there exist $R_0, \kappa > 0$ such that, for each $x\in\R^d$ with $|x|\geq R_0$,
    \[
        \langle b(x), x \rangle \leq -\kappa\, |x|\ .
    \]
\end{itemize}
Given the assumptions $(H)$ and $(AD)$,  the system \eqref{eq:system1} has a unique invariant measure $\mu$ (see Subsection \ref{subsec: Spectral properties unperturbed transop} for more details). Moreover, 
the dissipativity assumption implies that the flow associated to the vector field $b$ transports measures from outside the closed ball $\overline{B_{R_0}(0)}$ towards it.

\smallskip Let $x_0\in\R^d$ be a chosen point of the phase space and let  $g_{x_0}: \R^d\to\R\cup\{+\infty\}$ be such that
\begin{equation} \label{eq:distance1}
    g_{x_0}(x):=-\log d(x,x_0)\ ,
\end{equation}
with the convention that $g_{x_0}(x_0)=+\infty$, and $d$ is the Euclidean distance. This observable encodes the appearance of an extreme event for a given point $x_0$ in the phase space and measures how far we are from our chosen $x_0$ in a logarithmic scale. Throughout the paper $B_r(x_0)$ denotes the open Euclidean ball of radius $r>0$ centered at $x_0$. Let $h>0$ be a fixed time-step and let the sequence of discrete times in which we sample the process be denoted as $t_k:=kh$.

We now state the extreme value law described in the introduction.

\begin{restatable}{thm}{mymaintheorem}\label{teo:Teorema1}
    For fixed $h,\tau>0$ and $x_0\in\R^d$, let $X_t^x$ be the solution at time $t$ of the system \eqref{eq:system1}, whose drift satisfies assumptions $(H)$ and $(AD)$, and let $(u_n)_{n\in\N}\subset\R^+$ be such that
    \begin{equation}\label{eq:ConvergenzaMisuraPalle}
        n\, \mu (B_n) \to \tau\ ,
    \end{equation}
    as $n\to+\infty$, where $B_{n}$ is the ball $B_{e^{-u_n}}(x_0)$ and
    $t_{k}=kh$. Then
    \begin{equation}\label{eq:TesiTeorema1}
        \lim_{n\to+\infty}\mathbb{P} \otimes\mu\left(\left\{(\omega,x)\in \Omega\times\R^d : \max_{k=0,\dots,n-1}g_{x_0}(X_{t_k}^x(\omega))\leq u_{n}\right\}\right) = e^{-\tau}\ .
    \end{equation}
\end{restatable}

\begin{oss}\label{oss: Indipendenza da h}
    The right hand side of the asymptotic law found in \eqref{eq:TesiTeorema1} does not depend on $h$. Despite this independence, the limit as $h\to 0$ of the above result does not provide more information and should not be thought of as a way to get a continuous-time version of our result: the passage from discrete to continuous time might require another not obvious rescaling (see \cite[Appendix A]{EXT} and references therein for a thorough discussion of the matter). 
\end{oss}

An extreme value law of the same form was established in
\cite{EXT} under a quadratic dissipativity assumption.
Here we obtain this law for bounded drifts satisfying
\textup{(AD)}, using classical differentiable operator
perturbation theory.

\subsection{Organization of the paper} 

The plan of the paper is as follows: in Section \ref{sec: trans, reg} we set up the functional analytic framework necessary to study our problems and we perform some regularization estimates on the unperturbed and perturbed transfer operators which will be used in the sequel.  In Section \ref{sec: rare events}, we first recall the main abstract tools we use regarding classical operator perturbation theory and then we verify the necessary assumptions for the operators in order to  recover the laws of extreme events stated in Theorem \ref{teo:Teorema1}. Lastly, in Appendix \ref{app: StimaD derivata mista} we prove a mixed-derivative estimate of the heat kernel with constant diffusion. Together with the Aronson-type estimates provided by \cite[Theorem 1.2]{MPZ}, this result is necessary to establish the regularizing properties of the transfer operators.

\begin{notazz}
    We will denote by $L^p(\R^d)$ the space of measurable functions which are $p$-summable with respect to the Lebesgue measure on $\R^d$. We denote by $SM(\R^d)$ the set of Borel signed measures on $\R^d$ with bounded variation, i.e., $\|\eta\|=|\eta|(\R^d)<+\infty$ for each $\eta\in SM(\R^d)$.
\end{notazz}
    
\section{Transfer operators, Banach spaces and regularization lemmas}\label{sec: trans, reg}

In this section, starting from a general result on the transition probabilities of a system like \eqref{eq:system1} we define and study the basic properties of the transfer operators associated with such a system. In Section \ref{sec: Functional space} we introduce a certain weighted Banach space and we will see in Sections \ref{sec: Regularization properties unperturbed} and \ref{sec: Regularization properties perturbed} that the transfer operators when applied to these spaces have useful regularization properties.

\smallskip For the SDE considered in \eqref{eq:system1}, assuming the regularity assumption $(H)$ and the dissipative assumption $(AD)$, it is known (see, e.g. \cite{Friedman64, MPZRef26}) that the equation has the properties of strong existence  of the solutions and pathwise uniqueness. Moreover, Menozzi, Pesce and Zhang in \cite{MPZ} prove bounds on the transition probabilities for these systems in this setup. Such estimates imply that the transfer operator associated to the system has a regular kernel, hence regularizing properties. 

\smallskip Let $(\theta_t)_{t\geq 0}$ be the flow solving 
\[
    \begin{cases}
        \dot{\theta}_t(x)=b(\theta_t(x))\\
        \theta_0(x)=x
    \end{cases}\ ,
\]
where $b$ is the drift function as in \eqref{eq:system1}. For $\lambda>0$, $t>0$ and $x\in\R^d$, let us denote the Gaussian kernel by
\[
    g_\lambda(t,x)=t^{-\frac{d}{2}}\, e^{-\frac{\lambda}{t}|x|^2}\ .
\]

We require explicit bounds for the transition density $S_t$ of the system \eqref{eq:system1} and its derivatives. Since the drift $b$ is bounded, the flow does not provide additional information, as can be seen below. Therefore, the Aronson-type estimates on the kernel and its derivatives take a slightly different form compared to those in \cite[Theorem 7]{EXT}. The following result provides the specific bound, together with an estimate on the mixed derivative $|\nabla_x\nabla_y\, S_t(x,y)|$, which will allow us to establish the necessary Lasota-Yorke inequalities on an exponentially weighted Banach space. 

\begin{thm}\label{teo: MPZ Weak}
    Let us fix $T>0$. For each $0<t\leq T$ and $x\in \R^d$, let $X_t^x$ be the unique solution of \eqref{eq:system1} at time $t$. Then $X_t^x$ admits a transition density $S_t(x, y)$, for each $y\in\R^d$, which is of class $C^1$ in both variables. Moreover, for each $y\in\R^d$, the mixed derivative $\nabla_x\nabla_y\, S_t(x,y)$ exists and is continuous in $x$. Lastly, $S_t(x,y)$ satisfies the following estimates:

    \begin{itemize}
        \item[$(1)$] \textnormal{two sided density bounds}: there exist $\lambda_{0}\in (0,1]$ and $C_{0}\geq 1$ depending on $T,K,\|b\|_\infty$ and $d$ such that, for any $x,y\in \R^d$ and $0<t\leq T$,
        \[
            C_{0}^{-1}\, e^{-2\lambda_{0}^{-1}\|b\|_\infty^2 T}\, g_{\lambda_{0}^{-1}}(t/2,x-y) \leq S_t(x,y)\leq C_0\,  e^{\lambda_0\|b\|_\infty^2 T}\, g_{\lambda_{0}}(2t,x-y)\ .
        \]

        \item[$(2)$] \textnormal{gradient estimates}: there exist $\lambda_{1}\in (0,1]$ and $C_{1}\geq 1$ depending on $T,K,\|b\|_\infty$ and $d$ such that, for any $x,y\in \R^d$ and $0<t\leq T$
        \[
            |\nabla_{x}\, S_t(x,y)|\leq C_1\, e^{\lambda_1\|b\|_\infty^2 T}\, t^{-\frac{1}{2}}\, g_{\lambda_{1}}(2t,x-y)\ ,
        \]
        \[
            |\nabla_y\, S_t(x,y)|\leq C_1\, e^{\lambda_1 \|b\|_\infty^2 T}\, t^{-\frac{1}{2}}\, g_{\lambda_{1}}(2t,x-y)\ .
        \]

        \item[$(3)$] \textnormal{mixed derivative estimates}: there exist $\lambda_{2},C_{2}>0$ depending on $T,K,\|b\|_{\infty}$ and $d$ such that, for any $x,y\in \R^d$ and $0<t\leq T$
        \[
            |\nabla_{x}\nabla_y\, S_t(x,y)|\leq C_2\, e^{\lambda_2 \|b\|_\infty^2 T}\, t^{-1}g_{\lambda_{2}}(2t,x-y)\ .
        \]
    \end{itemize}

    \begin{myproof}
        The existence of the density and the estimates involving the flow $\theta_t(x)-y$ inside the Gaussian distribution follow directly from \cite[Theorem 7]{EXT} and Proposition \ref{prop: Stima Mixed Derivative} (see Appendix \ref{app: StimaD derivata mista}). Since the drift $b$ is bounded, we can remove the dependence on the flow by observing that:
        \begin{equation}\label{eq: No Additional Information}
            |x-y|-\|b\|_\infty\, t\leq |\theta_t(x)-y|\leq |x-y|+\|b\|_\infty\, t\ .
        \end{equation}
        The latter applies identically to the regularized flow $\theta_{t,0}^{(1)}(x)$ because the mollification does not increase $\|b\|_\infty$.
        
        For the upper bounds, we use the inequality $|\theta_t(x)-y|^2 \geq \frac{1}{2}|x-y|^2 - \|b\|_\infty^2 t^2$. Hence, for any $\lambda>0$:
        \[
            g_{\lambda}(t,\theta_t(x)-y) \leq t^{-\frac{d}{2}}\, e^{-\frac{\lambda}{t}\left(\frac{1}{2}\, |x-y|^2-\|b\|_\infty^2\, t^2\right)} = 2^{\frac{d}{2}}\, e^{\lambda\|b\|_\infty^2 t}\, (2t)^{-\frac{d}{2}}e^{-\frac{\lambda}{2t}|x-y|^2}\leq 2^{\frac{d}{2}}\, e^{\lambda\|b\|_\infty^2 T} g_\lambda(2t,x-y)\ .
        \]
        Applying this upper bound to the estimates from the aforementioned results yields the upper density bound, the gradient estimates and the mixed derivative estimates.

        For the lower density bound, we use the upper inequality from \eqref{eq: No Additional Information}, which implies $|\theta_t(x)-y|^2 \leq 2\, |x-y|^2 + 2\, \|b\|_\infty^2 t^2$. Therefore:
        \[
            g_{\lambda_0^{-1}}(t,\theta_t(x)-y) \geq t^{-\frac{d}{2}}\, e^{-\frac{\lambda_0^{-1}}{t}\left(2|x-y|^2 + 2\|b\|_\infty^2 t^2\right)}= 2^{-\frac{d}{2}}\, e^{-2\lambda_0^{-1}\|b\|_\infty^2 t}\, (t/2)^{-\frac{d}{2}}e^{-\frac{\lambda_0^{-1}}{\frac{t}{2}}|x-y|^2} \geq 2^{-\frac{d}{2}}\, e^{-2\lambda_0^{-1}\|b\|_\infty^2 T} g_{\lambda_0^{-1}}(t/2, x-y)\ .
        \]
        Applying this to the lower bound from \cite[Theorem 7]{EXT} completes the proof (note that all the constants $C_i$ absorb the factor $2^{\frac{d}{2}}$ and $2^{-\frac{d}{2}}$).
    \end{myproof}
\end{thm}

\subsection{The Kolmogorov operator and the transfer operator}\label{sec: Kolmogorov and Transfer unperturbed}

Here, as done in \cite[Section 3.1]{EXT}, we define the transfer operators associated to the evolution of system \eqref{eq:system1} and show some of the basic properties of these operators. The properties of the transition probabilities $S_t$ inherited by Theorem \ref{teo: MPZ Weak} will be used to define a Kolmogorov operator (composition operator) associated to our SDE.

\begin{defn}\label{def: koldef}
    The \textit{Kolmogorov operator} $P_{t}:L^\infty(\R^d)\to C^0(\R^d)$ associated to the system \eqref{eq:system1} and to the evolution time $t>0$ is defined  as follows. Let $\phi \in L^{\infty}(\R^d)$. Then for every $x \in \R^d$
    \[
    (P_{t}\phi)(x):=\mathbb{E}[\phi (X_t^x)]\ .
    \]
    In the literature, this is also known as ``\textit{stochastic Koopman operator}''. Given the transition probabilities $S_t$, it holds 
    \begin{equation}\label{eq: KolmogorovOperatorByMarkovKernel}
        (P_{t}\phi )(x)=\int_{\R^d}\phi(y)S_t(x,y)\, dy\ .
    \end{equation}
\end{defn}

We note that if $\nu $ is a Borel signed measure on $\R^{d}$, then
\begin{equation}\label{eq:DualityPerMisuraQualunque}
    \int_{\R^d} (P_t\phi )(x)\, d\nu (x)=\int_{\R^d} \left(\int_{\R^d}\phi(y)S_t(x,y)\, dy\right) d\nu (x)\ .
\end{equation}
Supposing that $\nu$ has a density $f\in L^1(\R^d)$ with respect to the Lebesgue measure, i.e.,  $d\nu = f(x)\, dx$, we can thus write 
\begin{equation}\label{eq:duality}
        \int_{\R^d} (P_{t}\phi )(x)\, d\nu(x) = \int_{\R^d} \phi (y) \left( \int_{\R^d} S_t(x,y)f(x)\, dx \right) dy\ .
\end{equation}
Therefore, it is possible to define the transfer operator $\cl_t$ as the dual of the Kolmogorov operator as follows.

\begin{defn}\label{def: TransOp}
    The \textit{unperturbed transfer operator} $\cl_t:L^1(\R^d)\to L^1(\R^d)$ associated to the system \eqref{eq:system1} and to the evolution time $t>0$ is defined as follows. Let $f\in L^1(\R^d)$. Then for almost every $y\in \R^d$ 
    \begin{equation}\label{eq:stoctransfer}
        (\cl_{t}\, f)(y):=\int_{\R^d} S_t(x,y)f(x)\, dx\ .
    \end{equation}
\end{defn}

\begin{oss}
    By \eqref{eq:duality}, we recover the duality relation between the Kolmogorov and the transfer operator whenever $d\nu = f(x) \, dx $, namely,
    \begin{equation}\label{eq:dualityfinale}
        \int_{\R^d} (P_{t}\phi )(x)f(x)\, dx=\int_{\R^d} \phi (y)(\cl_{t}\, f)(y)\, dy\ .
        \end{equation}
\end{oss}

The next lemma is a well-known property of the transfer operator.

\begin{lem}[Lemma 10, \cite{EXT}]\label{lem:TransOpWeakContraction}
    The unperturbed transfer operator $\cl_t$ preserves the integral and is a weak contraction with respect to the $L^1$-norm. Moreover, it is a Markov operator having kernel $S_t$.
\end{lem}

In the following lemma, we see how the presence of the noise provides a form of regularization for the transfer operator.

\begin{lem}\label{lem: C1Estimate}
    For every $t>0$, $\cl_t$ is a bounded linear operator from $L^1(\R^d)\to C^1_b(\R^d)$: in particular there exists a constant $C>0$ (depending on $t,K,\|b\|_\infty$ and $d$) such that, for each $f\in L^1(\R^d)$,
    \begin{equation}\label{eq:C1Estimate}
        \|\cl_t f\|_{C^1}\leq C\|f\|_{L^1}\ ,
    \end{equation}
    where $\|g\|_{C^1}=\|g\|_{\infty}+\|\nabla g\|_{\infty}$.
 
    \begin{myproof}
        For any $y\in\R^d$,
        \[  
            |(\cl_tf)(y)|\leq \int_{\R^d}S_t(x,y) |f(x)|\, dx\leq C_0e^{\lambda_0\|b\|_\infty^2t}(2t)^{-\frac{d}{2}}\|f\|_{L^1}\ ,
        \]
        where in the last estimate we bound the supremum norm of the function $S_t(\cdot\, , y)$ by Theorem \ref{teo: MPZ Weak} (the exponential in the Gaussian distribution $g_{\lambda_0}(2t,x-y)$ is bounded by $1$). 
        Moreover, by Theorem \ref{teo: MPZ Weak} the map $y\mapsto \nabla_y\, S_t(x,y)$ is continuous for each $x\in\R^d$ and uniformly bounded on $\R^d\times\R^d$. Since $f\in L^1(\R^d)$, the dominated convergence theorem justifies differentiation under integral sign and shows that $y\mapsto \nabla_y (\cl_t f)(y)$ is continuous on $\R^d$. Hence $\cl_{t}f\in C_b^1(\R^d)$, namely:
        \[  
            |\nabla_y(\cl_tf)(y)|\leq \int_{\R^d}|\nabla_y\, S_t(x,y)|\, |f(x)|\, dx\leq C_1e^{\lambda_1\|b\|_\infty^2t}t^{-\frac{1}{2}(d+1)}2^{-\frac{d}{2}}\|f\|_{L^1}\ .
        \]
        Therefore, 
        \[
            \|\cl_t f\|_{C^1}\leq (2t)^{-\frac{d}{2}}\left(C_0e^{\lambda_0\|b\|_\infty^2t}+C_1e^{\lambda_1\|b\|_\infty^2t}t^{-\frac{1}{2}}\right)\|f\|_{L^1}\ ,
        \]
        as wanted.
    \end{myproof}
\end{lem}

We note that Definition \ref{def: TransOp} applies to densities $f\in L^1(\R^d)$,  i.e., to measures which are absolutely continuous with respect to the Lebesgue measure. However, we can extend this definition to every Borel signed measure on $\R^d$ of bounded variation, even to those that are singular with respect to the Lebesgue measure. This should not seem strange:  it is suggested from \eqref{eq:DualityPerMisuraQualunque} and if we consider the system \eqref{eq:system1} in terms of the densities of the random variables involved, then its initial condition is indeed a Dirac measure. Hence, it is possible to give the following well-posed definition.

\begin{defn}\label{def: TransOpGeneralizzato}
    The \textit{unperturbed transfer operator} $\cl_t:SM(\R^d)\to L^1(\R^d)$ associated to the system \eqref{eq:system1} and to the evolution time $t>0$ is defined as follows. Let $\eta\in SM(\R^d)$ be a Borel signed measure on $\R^d$ with bounded variation. Then for almost every $y\in\R^d$
    \[
        (\cl_t\, \eta)(y):=\int_{\R^d}S_t(x,y)\, d\eta(x)\ .
    \]
\end{defn}

To prove the extreme event laws shown in Theorem \ref{teo:Teorema1} we consider the target set $B_n$ as a hole in the phase space and consider it as an open system. Then we will study it by means of the related transfer operators and their properties. We remark that we introduce the hole corresponding to the rare event \textit{before} the regularizing action of the unperturbed transfer operator, differently from what was done in \cite[Definition 12]{EXT}.

\begin{defn}\label{def: PertKol and PertTrans}
    Let $t>0$, $x_0\in \R^d $ and $(u_n)_{n\in\N}\subset\R$ such that $u_n\to +\infty$ for $n\to+\infty$. We denote by $B_n$ the ball $B_{e^{-u_n}}(x_0)$. 
    The \textit{perturbed Kolmogorov operator} is defined as follows. Let $\phi \in L^{\infty}(\R^d)$. Then for every $x \in \R^d$
    \begin{equation}\label{eq: KolmogorovperturbedOperatorByMarkovKernel} 
        (P_{t,n}\, \phi)(x):=\mathds{1}_{B_n^c}(x)(P_t\phi)(x)=\mathds{1}_{B_n^c}(x)\int_{\R^d}\phi(y)S_t(x,y)\, dy\ .
    \end{equation}
    The \textit{perturbed transfer operator} is defined as follows. Let  $f\in L^1(\R^d)$. Then for almost every $y\in\R^d$
    \begin{equation}\label{eq:TransPerturbed}
        (\cl_{t,n}\, f)(y):=\left(\cl_t\left(f\,  \mathds{1}_{B^c_n}\right)\right)(y)=\int_{\R^d} S_t(x,y)f(x)\mathds{1}_{B_n^c}(x)\, dx\ .
    \end{equation}
\end{defn}

\begin{oss}
    Using the definitions above, one can prove that the perturbed operators enjoy the following duality relation:
    \begin{equation} \label{eq:DualitaOperatoriPerturbati} 
            \int_{\R^d}(P_{t,n}\, \phi)(x)f(x)\, dx=\int_{\R^d}\phi(y)(\cl_{t,n}\, f)(y)\, dy\ .
    \end{equation}
\end{oss}

The iterates of the perturbed Kolmogorov operator inherit their properties from the following result.

\begin{lem}\label{lem: ComposizioneKolmogorov}
    For every $t,s>0$ and every $\phi\in L^\infty(\R^d)$,
    \[
        P_{t,n}(P_{s,n}\phi)(x) =\mathds{1}_{B_n^c}(x)\, \mathbb{E}\left[\mathds{1}_{B_n^c}(X_t^{x})\, \phi(X_{t+s}^x)\right]\ .
    \]
    
    \begin{myproof}
        The proof is very similar to the one of \cite[Lemma 13]{EXT}.
    \end{myproof}
\end{lem}

\begin{oss}
    In particular, given generic $t_1,t_2 \in \R^+$
    \[
        (P_{t_1+t_2,n}\phi)(x)\neq P_{t_1,n}(P_{t_2,n}\phi)(x)\ ,
    \] 
    namely $P_{t,n}$ is not a semigroup. However, notice that if $t_1<t_2$, then taking $t=t_1$ and $t+s=t_2$ above, by Lemma \ref{lem: ComposizioneKolmogorov} we get
    \[
        P_{t_1,n}(P_{t_2-t_1,n}\phi)(x)  =\mathds{1}_{B_n^c}(x)\, \mathbb{E}\left[\mathds{1}_{B_n^c}(X_{t_1}^x)\, \phi(X_{t_2}^x)\right]\  .
    \]
\end{oss}

Consequently, we obtain the following result for the $N$-th iterate of the Kolmogorov operator when the sequence of time steps is properly chosen. 

\begin{cor}\label{cor:EtoTwistedP}
    For every $x\in \R^d$, every sequence of times $0=t_0<t_1<\dots<t_N$ and every $\phi\in L^\infty(\R^d)$, we have
    \begin{equation}\label{eq:EtoTwistedP}
        \left(P_{t_1-t_0,n}\circ\cdots\circ P_{t_N-t_{N-1},n}\right)(\phi)(x) = \mathbb{E}\left[\mathds{1}_{B_n^c}(X_{t_{0}}^x) \cdots \mathds{1}_{B_n^c}(X_{t_{N-1}}^x)\, \phi(X_{t_N}^x)\right]\ .
    \end{equation}
\end{cor}

\subsection{Functional spaces}\label{sec: Functional space}

We now define suitable functional spaces on which the transfer operators introduced in the previous sections have a regularizing behavior and nice spectral properties. We construct spaces which are suitable for our non-compact environment, imposing a controlled behavior far away from the origin by using weight functions growing at infinity.

\smallskip Let $\alpha>0$. We consider the \textit{exponential weight function}
\[
    \rho_{\alpha,\exp}(x)=e^{\alpha\sqrt{1+|x|^2}}\ .
\]
We define $L^1_{\alpha,\exp}(\R^d)$ as the space of Lebesgue measurable functions $f:\R^d\to\R$ such that
\[
    \|f\|_{L^1_{\alpha,\exp}}:=\int_{\R^d}\rho_{\alpha,\exp}(x)|f(x)|\, dx<+\infty\ .
\]
We equip the intersection space $\Bexp(\R^d):=L^1_{\alpha,\exp}(\R^d)\cap C^1_b(\R^d)$ with the sum norm:
\[
    \|f\|_{\Bexp} = \|f\|_{L^1_{\alpha,\exp}} + \|f\|_{C^1}\ .
\]

\begin{oss}
    Note that $L^1_{\alpha,\exp}(\R^d)\subset L^1(\R^d)$ and if $f\in L^1_{\alpha,\exp}(\R^d)$, then $\|f\|_{L^1}\leq \|f\|_{L^1_{\alpha,\exp}}$ because $\rho_{\alpha,\exp}(x) \geq 1$. Moreover, equipped with this norm, $(\Bexp(\R^d), \|\cdot\|_{\Bexp})$ is a Banach space.
\end{oss}

\begin{prop}\label{prop: CompactEmbedding SpazioForte Weak}
    The embedding $\Bexp(\R^d)\hookrightarrow L^1(\R^d)$ is compact. 

    \begin{myproof}
        Since $\rho_{\alpha,\exp}(x)\geq e^\alpha\geq 1$ for all $x\in\R^d$, the inclusion is continuous. To show that it is compact, let $(f_n)_{n\in\N}$ be a bounded sequence in $\Bexp(\R^d)$. Because $\sup_{n\in\N}\|f_n\|_{C^1}\leq M$, $(f_n)_{n\in\N}$ is uniformly bounded and equicontinuous on $\R^d$, thus by Ascoli--Arzelà theorem, combined with a diagonal argument, there exists a subsequence $(f_{n_k})_{k\in\N}$ and a continuous function $f$ such that $f_{n_k}\to f$ uniformly on every compact subset of $\R^d$. By Fatou's lemma $\|f\|_{L^1_{\alpha,\exp}}\leq\liminf_{k\to+\infty}\|f_{n_k}\|_{L^1_{\alpha,\exp}}\leq M$, which implies that $f\in L_{\alpha,\exp}^1(\R^d)\subset L^1(\R^d)$ and $\|f_{n_k}-f\|_{L^1_{\alpha,\exp}}\leq 2M$ for each $k\in\N$. For any $R>0$, 
        \[
            \|f_{n_k}-f\|_{L^1}=\int_{|x|\leq R}|f_{n_k}(x)-f(x)|\, dx+\int_{|x|>R}|f_{n_k}(x)-f(x)|\, dx\leq \int_{|x|\leq R}|f_{n_k}(x)-f(x)|\, dx+2Me^{-\alpha\sqrt{1+R^2}}\ .
        \]
        Since $f_{n_k}\to f$ uniformly on the compact ball $\overline{B_R(0)}$, we have that $\limsup_{k\to+\infty}\|f_{n_k}-f\|_{L^1}\leq 2Me^{-\alpha\sqrt{1+R^2}}$. By arbitrariness of $R$, we get that $f_{n_k}$ converges to $f$ in the $L^1$-topology, as wanted.
    \end{myproof}
\end{prop}

\subsection{Regularizing properties for the unperturbed transfer operator}\label{sec: Regularization properties unperturbed}

In this section, we demonstrate how the dissipativity assumption and the noise in the SDE \eqref{eq:system1} induce a regularizing effect on the associated unperturbed transfer operator. With the ultimate goal of studying the extreme events of the discrete-time sampled process $X_{t_k}^x$ at times $t_k=kh$, we focus on the iterates $\cl_h^n$ of the transfer operator corresponding to the fixed time-step $t=h$.

\smallskip The main result of this section is a uniform Lasota-Yorke inequality for the unperturbed transfer operator $\cl_h$ acting on $\Bexp(\R^d)$ and $L^1(\R^d)$ as the strong and weak spaces, respectively. The starting point is Lemma \ref{lem:exponential_Lyapunov_estimate}, which establishes an exponential Lyapunov estimate under our current dissipativity condition. 

\begin{lem}\label{lem:exponential_Lyapunov_estimate}
    If $0<\alpha<2\kappa$, then there exist constants $\gamma,C_\alpha>0$ such that
    \begin{equation}\label{eq:exponential_Lyapunov_ineq}
        \mathcal{A}\, \rho_{\alpha,\exp}(x)\le-\gamma\,  \rho_{\alpha,\exp}(x)+C_\alpha 
    \end{equation}
    for every $x\in\mathbb R^d$, where $\mathcal A$ is the infinitesimal generator of the Markov semigroup associated with \eqref{eq:system1}, namely
    \begin{equation}\label{eq:generator_definition}
        \mathcal A\, \varphi= b\cdot \nabla\,\varphi+\frac{1}{2}\, \Delta\, \varphi\ ,
    \end{equation}
    for every smooth function $\varphi:\R^d\to\R$.
        
    \begin{myproof}
        Let $r(x)=(1+|x|^2)^{\frac{1}{2}}$. Direct computation yields the gradient and Laplacian of $r$:
        \[
            \nabla\, r(x) = \frac{x}{r(x)}\ , \qquad \Delta\, r(x) = \frac{d}{r(x)} - \frac{|x|^2}{r(x)^3}\ .
        \]
        In particular, we have the bounds $|\nabla\, r(x)|^2 \le 1$ and $\Delta\, r(x) \le \frac{d}{r(x)}$. 
        
        Applying the chain rule to the weight $\rho_{\alpha,\exp}(x)=e^{\alpha\, r(x)}$, we find:
        \[
            \nabla\, \rho_{\alpha,\exp}(x) = \alpha\, \rho_{\alpha,\exp}(x)\nabla\, r(x)\ ,\qquad
            \Delta\, \rho_{\alpha,\exp}(x)= \alpha\, \rho_{\alpha,\exp}(x)\Delta\, r(x) + \alpha^2\, \rho_{\alpha,\exp}(x)|\nabla\, r(x)|^2\ .
        \]
        Evaluating the infinitesimal generator $\mathcal{A}$ on $\rho_{\alpha,\exp}$ and dividing by the strictly positive weight itself yields:
        \begin{equation}\label{eq:AW_over_W}
            \frac{\mathcal{A}\, \rho_{\alpha,\exp}(x)}{\rho_{\alpha,\exp}(x)} = \alpha\, \langle b(x),\nabla\, r(x)\rangle + \frac{\alpha}{2}\, \Delta\, r(x) + \frac{\alpha^2}{2}\, |\nabla\, r(x)|^2\ .
        \end{equation}
        We now apply the dissipativity assumption $(AD)$: for every $x\in\R^d$ such that $|x|\ge R_0$, we have
        \[
            \langle b(x),\nabla\, r(x)\rangle = \frac{\langle b(x),x\rangle}{r(x)} \le -\frac{\kappa\, |x|}{r(x)}\ .
        \]
        Substituting this and the differential bounds into \eqref{eq:AW_over_W}, we obtain, for all $|x|\ge R_0$,
        \begin{equation}\label{eq:AW_over_W_bound}
            \frac{\mathcal{A}\, \rho_{\alpha,\exp}(x)}{\rho_{\alpha,\exp}(x)} \le - \frac{\alpha\, \kappa\, |x|}{r(x)} + \frac{\alpha\, d}{2\, r(x)} + \frac{\alpha^2}{2}\ .
        \end{equation}
        As $|x|\to+\infty$, the ratio $\frac{|x|}{r(x)} \to 1$ while $\frac{1}{r(x)} \to 0$. Therefore,
        \[
            \limsup_{|x|\to+\infty} \frac{\mathcal{A}\, \rho_{\alpha,\exp}(x)}{\rho_{\alpha,\exp}(x)} \le -\alpha\kappa + \frac{\alpha^2}{2}\ .
        \]
        Because $0<\alpha<2\kappa$ by hypothesis, this limit is strictly negative. Consequently, there exists a constant $\gamma>0$ and a sufficiently large radius $R_1\ge R_0$ such that $\mathcal{A}\, \rho_{\alpha,\exp}(x) \le -\gamma\, \rho_{\alpha,\exp}(x)$ for each $|x| \ge R_1$.

        Finally, it remains to control the bound on the compact region $\overline{B_{R_1}(0)}$. Because the drift $b$ is continuous and $\rho_{\alpha,\exp}$ is smooth, the function  $x\mapsto \mathcal{A}\, \rho_{\alpha,\exp}(x) + \gamma\, \rho_{\alpha,\exp}(x)$ is continuous on the compact ball $\overline{B}_{R_1}(0)$, thus it is bounded from above by some constant $C_\alpha > 0$. Combining this bound with the decay condition outside the ball ensures that \eqref{eq:exponential_Lyapunov_ineq} holds globally on $\R^d$.
    \end{myproof}
\end{lem}

The proof of the following result adapts the strategy of \cite[Lemma 18]{EXT} to our current functional setting, relying crucially on the exponential Lyapunov estimate established in the previous lemma.

\begin{lem}\label{lem: LSY L12 Weak}
    For every $h>0$ and $0<\alpha<2\kappa$, there exist constants $A,B>0$ and $\lambda\in(0,1)$ (depending on $h$) such that, for each $f\in L^1_{\alpha,\exp}(\R^d)$ and $n\in\N$,
    \begin{equation}\label{eq:LSYforL12 Weak}
        \|\cl_h^n f\|_{L^1_{\alpha,\exp}}\leq A\lambda^n\|f\|_{L^1_{\alpha,\exp}}+B\|f\|_{L^1}\ .
    \end{equation} 

    \begin{myproof}
        Let $L_{\mathrm{dens}}^1(\R^d)$ denote the set of probability density functions. By standard decomposition of $f\in L^1_{\alpha,\exp}(\R^d)$ into its positive and negative parts, it suffices to prove that there exist $C, D>0$ such that 
        \begin{equation}\label{eq:ineq2 lemma 18 Weak}
            \| \cl_h^n p\| _{L^1_{\alpha,\exp}}\leq C\lambda^{n}\| p\| _{L^1_{\alpha,\exp}}+D  
        \end{equation}
        for every $p\in L_{\text{dens}}^1(\R^d)\cap L^1_{\alpha,\exp}(\R^d)$. 
        
        By the Fubini--Tonelli theorem, we can rewrite \eqref{eq:ineq2 lemma 18 Weak} in terms of expectations. Specifically, the latter is equivalent to showing:
        \begin{equation}\label{eq:ineq3 lemma 18 Weak}
            \mathbb{E}\left[\rho_{\alpha,\exp}(X_{nh})\right] \leq C\lambda ^{n}\mathbb{E}\left[\rho_{\alpha,\exp}(X_0)\right] +D\ ,
        \end{equation}
        where $X_0$ is a random initial condition with density $p$ and $X_t$ is the solution of the Cauchy problem \eqref{eq:system1}.
        
        Since $\rho_{\alpha,\exp}$ is smooth, we can apply It\^{o}'s formula to $(t, x) \mapsto e^{\gamma t}\rho_{\alpha,\exp}(x)$. Given $R>0$, let $\tau_R:=\inf\{t\geq 0 : |X_t|\geq R\}$ 
        (with the convention that $\inf\emptyset=+\infty$). For the stopped process, we obtain:
        \begin{align*}
            e^{\gamma (\min\{t,\tau_R\})}\rho_{\alpha,\exp}(X_{\min\{t,\tau_R\}}) &= \rho_{\alpha,\exp}(X_0) + \int_0^{\min\{t,\tau_R\}} e^{\gamma s} \left(\gamma \rho_{\alpha,\exp}(X_s) + \mathcal{A}\, \rho_{\alpha,\exp}(X_s)\right) ds + \int_0^{\min\{t,\tau_R\}} e^{\gamma s} \nabla \rho_{\alpha,\exp}(X_s)\cdot dW_s\ .
        \end{align*}
        Notice that for any $s\in (0,\min\{t,\tau_R\})$, we necessarily have $\tau_R>0$ and $|X_s|<R$ (while on $\{\tau_R=0\}$ the stochastic integral vanishes identically). Therefore, the integrand $\mathds{1}_{(0,\tau_R]}(s)\, e^{\gamma s}\, \nabla\rho_{\alpha,\exp}(X_s)$ is uniformly bounded on $(0,t]$, so the Itô integral is a true martingale and its expectation is zero. Recalling Lemma \ref{lem:exponential_Lyapunov_estimate}, we have $\mathcal{A}\, \rho_{\alpha,\exp}(x) \leq -\gamma\, \rho_{\alpha,\exp}(x)+C_\alpha$, from which we deduce that $\gamma \rho_{\alpha,\exp}(x) + \mathcal{A}\, \rho_{\alpha,\exp}(x) \leq C_\alpha$. Taking the expected value yields:
        \begin{align*}
            \mathbb{E}\left[e^{\gamma (\min\{t,\tau_R\})}\rho_{\alpha,\exp}(X_{\min\{t,\tau_R\}})\right] &\leq \mathbb{E}\left[\rho_{\alpha,\exp}(X_0)\right] + \mathbb{E}\left[\int_0^{\min\{t,\tau_R\}} C_\alpha e^{\gamma s}\, ds\right] \\
            &\leq \mathbb{E}\left[\rho_{\alpha,\exp}(X_0)\right] + \int_0^t C_\alpha e^{\gamma s}\, ds \\
            &= \mathbb{E}\left[\rho_{\alpha,\exp}(X_0)\right] + \frac{C_\alpha}{\gamma}(e^{\gamma t}-1)\ .
        \end{align*}
        Since this bound is independent of $R$, we can apply Fatou's lemma as $R \to +\infty$. Note that because the solution $X_t$ exists globally, $\lim_{R\to+\infty}\tau_R=+\infty$ almost surely:
        \[
            e^{\gamma t}\mathbb{E}\left[\rho_{\alpha,\exp}(X_t)\right] = \mathbb{E}\left[\liminf_{R\to+\infty} e^{\gamma (\min\{t,\tau_R\})}\rho_{\alpha,\exp}(X_{\min\{t,\tau_R\}})\right] \leq \liminf_{R\to+\infty} \mathbb{E}\left[e^{\gamma (\min\{t,\tau_R\})}\rho_{\alpha,\exp}(X_{\min\{t,\tau_R\}})\right]\ .
        \]
        Therefore,
        \[
            e^{\gamma t}\mathbb{E}\left[\rho_{\alpha,\exp}(X_t)\right] \leq \mathbb{E}\left[\rho_{\alpha,\exp}(X_0)\right] + \frac{C_\alpha}{\gamma}e^{\gamma t}\ .
        \]
        Dividing both sides by $e^{\gamma t}$ and evaluating at $t=nh$, we get \eqref{eq:ineq3 lemma 18 Weak}, as wanted.
    \end{myproof}
\end{lem}

Now we are able to prove the following Lasota-Yorke inequality for $\cl_h$ with the spaces $\Bexp(\R^d)$ and $L^1(\R^d)$.

\begin{prop}\label{prop: LYITransWeak}
    For every $h>0$ and $0<\alpha<2\kappa$, there exist constants $\tilde{A},\tilde{B}>0$ and $\lambda\in (0,1)$ (depending on $h$) such that for each $f\in \Bexp(\R^d)$ and $n\geq 1$
    \begin{equation}\label{eq:LSYTransWeak}
        \|\cl_h^nf\|_{\Bexp}\leq \tilde{A}\lambda^n\|f\|_{\Bexp}+\tilde{B}\|f\|_{L^1}\ .
    \end{equation}

    \begin{myproof}
        We have 
        \begin{align*}
            \|\cl_h^nf\|_{\Bexp}&\leq \left(A\lambda^n\|f\|_{L^1_{\alpha,\exp}}+B\|f\|_{L^1}\right)+\|\cl_h(\cl_h^{n-1}f)\|_{C^1}\leq \left(A\lambda^n\|f\|_{\Bexp}+B\|f\|_{L^1}\right)+C\|\cl_h^{n-1}f\|_{L^1}\\
            &\leq \left(A\lambda^n\|f\|_{\Bexp}+B\|f\|_{L^1}\right)+C\|f\|_{L^1}=A\lambda^n\|f\|_{\Bexp}+(B+C)\|f\|_{L^1}\ ,
        \end{align*}
         where we applied Lemma \ref{lem: LSY L12 Weak} to bound $\|\cl_h^n f\|_{L^1_{\alpha,\exp}}$ and Lemma \ref{lem: C1Estimate} to estimate $\|\cl_h^n f\|_{C^1}=\|\cl_h(\cl_h^{n-1}f)\|_{C^1}$, which gives \eqref{eq:LSYTransWeak} with $\tilde{A}=A$ and $\tilde{B}=B+C$.
    \end{myproof}
\end{prop}

\subsection{Regularization properties for the perturbed transfer operators}\label{sec: Regularization properties perturbed}

In this section, we extend the regularization properties to the perturbed transfer operators. We prove a uniform Lasota-Yorke inequality for $\cl_{h,n}$ acting on the strong space $\Bexp(\R^d)$ and the weak space $L^1(\R^d)$. As before, the time evolution is fixed to $h>0$.

\begin{lem}\label{lem: LYIExpPerturbedTrans}
    For every $h>0$ and $0<\alpha<2\kappa$, there exist constants $A,B>0$ and $\lambda\in (0,1)$ (depending on $h$) such that for each $f\in L^1_{\alpha,\exp}(\R^d)$ and $n,m\geq 1$:
    \begin{equation}\label{eq:ineq1LYIExpPerturbedTrans}
        \|\cl_{h,n}^m\, f\|_{L^1_{\alpha,\exp}}\leq A\lambda^m\|f\|_{L^1_{\alpha,\exp}}+B\|f\|_{L^1}\ .
    \end{equation}

    \begin{myproof}
        We first assume $f\geq 0$. By definition, $\cl_{h,n}\, f=\cl_h(f\, \mathds{1}_{B_n^c})$, meaning $\cl_h\, f=\cl_{h,n}\, f+\cl_h(f\, \mathds{1}_{B_n})$. For fixed $n\geq 1$, let $g_{m}=\cl_h^m\, f- \cl_{h,n}^m\, f$. Because $\cl_h$ is a positive operator and $f\geq 0$, by induction it follows that $g_m\geq 0$. Consequently, $0 \leq \cl_{h,n}^{m}\, f \leq \cl_{h,n}^{m}\, f+g_{m} = \cl_h^m\, f$. Taking the $L^1_{\alpha,\exp}$-norm preserves this inequality: $\|\cl_{h,n}^{m}\, f\|_{L^1_{\alpha,\exp}} \leq \|\cl_h^m\, f\|_{L^1_{\alpha,\exp}}$.
        Applying Lemma \ref{lem: LSY L12 Weak} to bound the unperturbed operator yields:
        \begin{equation}\label{eq:LSYExp}
             \|\cl_{h,n}^{m}\, f\|_{L^1_{\alpha,\exp}} \leq A\lambda^m\|f\|_{L^1_{\alpha,\exp}}+B\|f\|_{L^1}\ .
        \end{equation}
        
        For a general $f \in L^1_{\alpha,\exp}(\R^d)$, we decompose $f=f^+-f^-$ into its positive and negative parts. By the linearity of $\cl_{h,n}^m$ and the triangle inequality, applying \eqref{eq:LSYExp} to $f^+$ and $f^-$ separately yields the desired bound, using that $f^+$ and $f^-$ have disjoint supports.
    \end{myproof}
\end{lem}

\begin{prop}\label{prop: LYIPerturbedTransSpazioForteExp}
    For every $h>0$ and $0<\alpha<2\kappa$, there exist constants $\tilde{A},\tilde{B}>0$ and $\lambda\in (0,1)$ (depending on $h$) such that for each $f\in \Bexp(\R^d)$ and $n,m\geq 1$:
    \begin{equation}\label{eq:LYIPerturbedTransSpazioForteExp}
        \|\cl_{h,n}^m\, f\|_{\Bexp}\leq \tilde{A}\lambda^m\|f\|_{\Bexp}+\tilde{B}\|f\|_{L^1}\ .
    \end{equation}

    \begin{myproof}
        By Lemma \ref{lem: LYIExpPerturbedTrans} and the bound $\|f\|_{L^1_{\alpha,\exp}} \leq \|f\|_{\Bexp}$, we immediately have:
        \begin{equation}\label{eq:ineq1LemmaLYIPeturbedTransExp}
            \|\cl_{h,n}^m\, f\|_{L^1_{\alpha,\exp}} \leq A\lambda^m\|f\|_{\Bexp}+B\|f\|_{L^1}\ .
        \end{equation}
        
        To control the $C^1$-component of the $\Bexp$-norm, we apply Lemma \ref{lem: C1Estimate} to $(\cl_{h,n}^{m-1}\, f) \mathds{1}_{B_n^c}\in L^1(\R^d)$ and get:
        \begin{equation}\label{eq: ineq2LemmaLYIPerturbedTrans}
            \|\cl_{h,n}^m\, f\|_{C^1} = \left\|\cl_h \left((\cl_{h,n}^{m-1}\, f) \mathds{1}_{B_n^c}\right)\right\|_{C^1} \leq C\left\|(\cl_{h,n}^{m-1}\, f) \mathds{1}_{B_n^c}\right\|_{L^1} \leq C\|\cl_{h,n}^{m-1}\, f\|_{L^1}\leq C\|f\|_{L^1} ,
        \end{equation}
        where the last inequality follows by induction from the fact that $\cl_{h,n}$ is an $L^1$-contraction.
        Adding the bounds \eqref{eq:ineq1LemmaLYIPeturbedTransExp} and \eqref{eq: ineq2LemmaLYIPerturbedTrans} establishes the Lasota-Yorke inequality of the statement with $\tilde{A}=A$ and $\tilde{B}=B+C$.
    \end{myproof}
\end{prop}

\section{Rare events via transfer operator}\label{sec: rare events}

For completeness, we recall the statement of two results due to Hennion--Hervé \cite{HennionHerve} which are fundamental to our construction.

\smallskip Let $(\B,\|\cdot\|)$ be a Banach space, $\B^\ast$ be its dual and $\Ll(\B)$ be the space of all bounded linear operators on $\B$. Let $\mathcal{V}$ denote the collection of all open intervals in $\R$ centered at zero and define $(v\otimes \varphi)f=\langle\varphi,f\rangle\, v$ for each $v,f\in\B$ and $\varphi\in\B^\ast$. Whenever spectral properties are considered in Subsections \ref{subsec: Spectral properties unperturbed transop} and \ref{subsec: Final proof}, we work on the standard complexification of the real Banach space $\Bexp(\R^d)$, i.e., we allow function to be complex-valued, extending their transfer operators by complex linearity. All estimates in Section \ref{sec: trans, reg} and Subsections \ref{subsec: Riparametrizzazione} and \ref{subsec: Derivate Transop} are carried out on the real space but extend automatically in this setting. 

\begin{defn}\label{defn: Quasi Compatto}
    We say that a bounded linear operator $Q:\B\to\B$ is \textit{quasi-compact} if there exist two $Q$-invariant closed subspaces such that $\B=F\oplus H$, where 
    \begin{itemize}
        \item $r(Q_{|H})<r(Q)$;
        \item $\dim\, F<+\infty$ and each eigenvalue of $Q_{|F}$ has modulus equal to $r(Q)$.
    \end{itemize}
\end{defn}

The following result is a generalization of the Ionescu--Tulcea--Marinescu theorem (see \cite{marinescu1950theorie}) which ensures quasi-compactness.

\begin{thm}[Theorem II.5, \cite{HennionHerve}]\label{teo: QC}
    Let $\|\cdot\|'$ be a continuous semi-norm on $\B$ and $Q:\B\to\B$ be a bounded linear operator such that
    \begin{itemize}
        \item [$(a)$] $Q(\{f\in\B : \|f\|\leq 1\})$ is conditionally compact in $(\B,\|\cdot\|')$, i.e., for any sequence $(f_n)_{n\geq 1}\subset\B$ such that $\|f_n\|\leq 1$ for every $n\geq 1$, there exist a subsequence $(f_{n_k})_{k\geq 1}\subset\B$ and $h\in\B$ such that 
        \[
        \lim_{k\to+\infty}\|Qf_{n_k}-h\|'=0\ ;
        \]
        \item [$(b)$] there exists a constant $M\geq 0$ such that 
        \[
        \|Qf\|'\leq M\|f\|'
        \]
        for every $f\in\B$;
        \item [$(c)$] there exist $k\in\N^\ast$ and $r,R\in\R^+$ such that
        \[
            r<r(Q)\ , \qquad \|Q^kf\|\leq r^k\|f\|+R\|f\|'
        \]
        for every $f\in\B$.
    \end{itemize}
    Then $Q$ is quasi-compact.
\end{thm}

As it will be seen below, when the peripheral eigenvalues of a quasi-compact operator $Q$ of diagonal type (i.e., the restriction $Q_{|F}$ is diagonalizable) are simple, smooth perturbations of this operator are still quasi-compact of diagonal type. However, the number of peripheral eigenvalues is not necessarily preserved. For this reason, it is convenient to introduce a new definition.

\begin{defn}
    Let $s\in\N^\ast$. We say that a bounded linear operator $Q:\B\to\B$ has \textit{$s$ dominating simple eigenvalues} if there exist two $Q$-invariant closed subspaces
    such that $\B=F\oplus H$, where
    \begin{itemize}
        \item $\dim\, F=s$ and $Q_{|F}$ has $s$ simple eigenvalues $\lambda_k$;
        \item $r(Q_{|H})<\min\{|\lambda_k| : k=1,\dots,s\}$.
    \end{itemize} 
\end{defn}

The following is a perturbation result and its corollary, useful for the subsequent construction.

\begin{thm}[Theorem III.8, \cite{HennionHerve}]\label{teo: TeoPerturbativo}
    Let $m \in \mathbb{N}^*$, $I_0\in\mathcal{V}$ centered at zero and $(Q(t))_{t \in I_0}$ be a collection of bounded linear operators on $\B$ such that 
    \begin{description}
        \item [${(H1)}$] $Q(\cdot) \in \mathcal{C}^m(I_0, \Ll(\mathcal{B}))$;
        \item [${(H2)}$] $Q(0)$ has $s$ dominating simple eigenvalues with $r(Q(0))=1$.
    \end{description}
    Then there exists $I_1\in\mathcal{V}$ with $I_1 \subset I_0$ such that $Q(t)$ has $s$ dominating simple eigenvalues for every $t \in I_1$. More precisely, there exist $\eta_1,\eta_2>0$, $c_1\geq 0$ and functions $\lambda_k(\cdot)\in C^m(I_1,\C)$, $v_k(\cdot)\in C^m(I_1,\B)$, $\varphi_k(\cdot)\in C^m(I_1,\B^\ast)$ for $k=1,\dots,s$ and $N(\cdot)\in C^m(I_1,\Ll(\B))$ such that, for every $t\in I_1$ and $k,l=1,\dots,s$,
    \begin{itemize}
        \item[${(i)}$] $Q(t)v_k(t) = \lambda_k(t)v_k(t)$ and $Q(t)^*\varphi_k(t) = \lambda_k(t)\varphi_k(t)$;
        \item[${(ii)}$] $\min\{|\lambda_k(t)| : k=1,\dots,s\}\geq 1-\eta_1$;
        \item[${(iii)}$] $\langle\varphi_k(t), v_l(t)\rangle = \delta_{k,l}$;
        \item[${(iv)}$] $Q^n(t) = \sum_{k=1}^s \lambda_k(t)^n (v_k(t) \otimes \varphi_k(t)) + N(t)^n$ for every $n\geq 1$ with
        \[
        \max\left\{\left\|\frac{d^l}{dt^l}N(t)^n\right\| : l=0,\dots,m\right\}\leq c_1(1-\eta_1-\eta_2)^n\ .
        \]
    \end{itemize}
\end{thm}

\begin{cor}[Corollary III.11, \cite{HennionHerve}]\label{cor:derivative}
    Under the setting of Theorem \ref{teo: TeoPerturbativo}, with $m \ge 1$, we have
    \[
        \lambda_1'(0) =\langle\varphi_1(0), Q(\cdot)'(0)v_1(0)\rangle\ .
    \]
\end{cor}

\subsection{Re-parameterization by target volume}\label{subsec: Riparametrizzazione}

Since Theorem \ref{teo: TeoPerturbativo} requires a family of operators depending on a parameter which varies in an open neighborhood of zero, we have to re-parameterize the family $\cl_h$ and $(\cl_{h,n})_{n\geq 1}$. Firstly, let us associate the unperturbed operator with $\cl_h$, i.e., we set 
\begin{equation}\label{eq: Riparametrizzazione unperturbed}
    \tilde{\cl}_{h,0}:=\cl_h\ .
\end{equation}
Subsequently, since the perturbed transfer operators differ from one another by the radius of the ball over which we integrate, we define a parametrization that interpolates this countable family of radii. Let $x_0\in\R^d$ be a chosen point in the phase space, $\omega_{d-1}$ be the surface area of the unit sphere in $\R^d$, $f_0\in\Bexp(\R^d)$ be the unique strictly positive\footnote{It turns out that $f_0(y)=(\tilde{\cl}_{h,0}f_0)(y)=\int_{\R^d}S_h(x,y) f_0(x)\, dx$ with $f_0\geq 0,\int_{\R^d}f_0(x)\, dx=1$ and $S_h(x,y)>0$ for each $x,y\in\R^d$ by Theorem \ref{teo: MPZ Weak}, which implies $f_0(y)>0$ for each $y\in\R^d$.} invariant density of the unique stationary measure $\mu$ of system \eqref{eq:system1} (see Proposition \ref{prop: AutovaloreUnperturbedTransOp} below) and $\tau>0$ as in \eqref{eq:ConvergenzaMisuraPalle}. Then, denoting by $I_0\subset\R$ a bounded open interval  centered at zero, we define
\begin{equation}\label{eq:ChosenRadiusFunction}
    \funcDef{r}{I_0\smallsetminus\{0\}}{\R^+}
    {\eps}{r(\eps):=\left(\dfrac{\tau\, d}{f_0(x_0)\, \omega_{d-1}}\, |\eps|\right)^{\frac{1}{d}}}\ ,
\end{equation}
and we extend it by continuity to all of $I_0$ by setting\footnote{This is consistent with the fact that $\tilde{\cl}_{h,0}$ integrates over $\R^d$, which is the same as integrating over $\mathbb{R}^d\smallsetminus\{x_0\}$, given that the singleton $\{x_0\}\subset\R^d$ has zero Lebesgue measure.} $r(0):=0$.

\begin{oss}
    Let us denote by $B_{r(\eps)}$ the ball $B_{r(\eps)}(x_0)$. Since $f_0$ is continuous at $x_0$, the Lebesgue differentiation theorem implies that
    \begin{equation}\label{eq: LDT applied to ball measure}
        \mu(B_{r(\eps)})=f_0(x_0)\, m(B_{r(\eps)})+o(r(\eps)^d)=f_0(x_0)\, \frac{\omega_{d-1}}{d}\left(\dfrac{\tau\, d}{f_0(x_0)\, \omega_{d-1}}\, |\eps|\right)+o(|\eps|)=\tau\, |\eps|+o(|\eps|)\ . 
    \end{equation}
     Dividing by $|\eps|$ and taking the limit as $\eps \to 0$, we conclude that
     \begin{equation}\label{eq:VelocitaConvergenzaPalleEps}
        \frac{\mu(B_{r(\eps)})}{|\eps|}\to\tau\ .
    \end{equation}
\end{oss}

Let us set $M_0=\Id$ and, for each $\eps\in I_0\smallsetminus\{0\}$,
\begin{equation}\label{eq:OperatoreBucoContinuo}
    \funcDef{M_\eps}{L^1(\R^d)}{L^1(\R^d)}{f}{M_\eps\, f:= f\, \mathds{1}_{B_{r(\eps)}^c}}\ .
\end{equation}
In this way, we can extend the definition of the discrete family of perturbed transfer operators $(\cl_{h,n})_{n\geq 1}$ into a continuous family, parameterized by the positive radius of the target ball whose complement we integrate over. Namely, if we denote $I_0^-=I_0\cap\R^-$ and $I_0^+=I_0\cap\R^+$, then we can define the family $(\tilde{\cl}_{h,\eps})_{\eps\in I_0}$ as follows:
\begin{equation}\label{eq:DefTransOpContinuo}
    \tilde{\cl}_{h,\eps}:=
        \begin{cases}
            2\cl_h-\cl_{h,r(-\eps)} & \eps\in I_0^-\\
            \cl_h & \eps=0\\
            \cl_{h,r(\eps)} & \eps\in I_0^+
        \end{cases}\ ,
\end{equation}
where, for each $\xi\in I_0\smallsetminus\{0\}$, 
\[
    (\cl_{h,r(\xi)}\, f)(y):=(\cl_h(M_\xi\, f))(y)=\int_{\R^d}S_h(x,y)(M_\xi\, f)(x)\, dx=\int_{B_{r(\xi)}^c}S_h(x,y)f(x)\, dx\ .
\]
Therefore,
\begin{equation}\label{eq: Def tranop continuo esplicito}
    (\tilde{\cl}_{h,\eps}\, f)(y)=\begin{cases}
        \displaystyle\int_{\R^d} S_h(x,y)f(x)\, dx+\displaystyle\int_{B_{r(\eps)}}S_h(x,y)f(x)\, dx & \eps\in I_0^-\\[.3cm]
        \displaystyle\int_{\R^d}S_h(x,y)f(x)\, dx\ & \eps=0\\[.3cm]
        \displaystyle\int_{B_{r(\eps)}^c}S_h(x,y)f(x)\, dx & \eps\in I_0^+
    \end{cases}\ .
\end{equation}

\begin{oss}\label{oss: Passaggio da famiglia continua a discreta}
    It will be useful to interpret the family $(\tilde{\cl}_{h,\eps})_{\eps\in I_0}$ as the following composition:
    \begin{equation}\label{eq: DefinizionTransFinaleComeComposizione}
        \tilde{\cl}_{h,\eps}=\begin{cases}
            \tilde{\cl}_{h,0}\circ (2M_0-M_{-\eps}) & \eps\in I_0^-\\
            \tilde{\cl}_{h,0}\circ M_0 & \eps=0\\
            \tilde{\cl}_{h,0}\circ M_\eps & \eps\in I_0^+    
        \end{cases}\ .
    \end{equation}
    Moreover, the discrete family $(\cl_{h,n})_{n\geq1}$ is recovered from the family $(\tilde{\cl}_{h,\eps})_{\eps\in I_0\smallsetminus\{0\}}$ at the values
    \begin{equation}\label{eq: Eps_n}
        \eps_n=\frac{f_0(x_0)\, \omega_{d-1}}{\tau\, d}\, e^{-d\, u_n}\ .
    \end{equation}
    Since $u_n\to+\infty$, we have $\eps_n\to 0$ as $n\to+\infty$, hence this sequence of points eventually falls into $I_0$ and it is such that $r(\eps_n)=e^{-u_n}$ for every $n\geq 1$, as wanted.
\end{oss}

\subsection{Verifying assumption \texorpdfstring{$(H1)$}{H1}}

Here our goal is to verify that the map $I_0\ni\eps \mapsto \tilde{\cl}_{h,\eps}\in\Ll(\Bexp(\R^d))$ is of class $C^1$, checking off the hypothesis $(H1)$ of Theorem \ref{teo: TeoPerturbativo}. 

\subsubsection*{Kantorovich--Rubinstein derivative of the intermediate operator \texorpdfstring{$M_\eps$}{M-epsilon}}\label{subsec: Derivata operatore intermedio}

Before doing this, we must understand how the target hole itself changes as the parameter $\eps$ varies. We compute the derivative of the intermediate hole operator $M_\eps$, with respect to $\eps$ in the Kantorovich--Rubinstein metric, i.e., we compute an operator $M'(\eps)$, at $\eps=0$ and at $\eps\in I_0^+$, such that, respectively,
\begin{equation}\label{eq: DerivataBucoContinuo}
    \lim_{\xi\to 0^+}\left\|\frac{(M_{\xi}-M_0)f}{\xi}-M'(0)f\right\|_{KR}=0\ ,\qquad \lim_{\xi\to 0}\left\|\frac{(M_{\eps+\xi}-M_\eps)f}{\xi}-M'(\eps)f\right\|_{KR}=0\ , 
\end{equation}
where we recall that the Kantorovich--Rubinstein norm on $\R^d$ is defined as follows: let $\mu\in SM(\R^d)$ be a Borel signed measure on $\R^d$ with bounded variation and $g:\R^d\to\R$. Then
\[
    \|\mu\|_{KR}:=\sup\left\{\int_{\R^d}g\, d\mu : \|g\|_\infty\leq 1,\,  \Lip(g)\leq 1\right\}\ ,
\]
where $\Lip(g)$ denotes the best Lipschitz constant of $g$. Thus, if $\mu,\nu\in SM(\R^d)$, their Kantorovich--Rubinstein distance is given by:
\begin{equation}\label{eq: def Kantorovich--Rubinstein distance}
    \|\mu-\nu\|_{KR}=\sup\left\{\int_{\R^d}g\, d(\mu-\nu): \|g\|_\infty\leq 1,\, \Lip(g)\leq 1\right\}\ . 
\end{equation}
Using the coarea formula and the properties of the Hausdorff measure, we establish the convergence of \eqref{eq: DerivataBucoContinuo} both at $\eps=0$ and at $\eps\in I_0^+$.

\begin{prop}\label{prop:ConvergenzaWassBuco}
    Let $x_0\in\R^d$, $(M_\eps)_{\eps\in I_0^+}$ be the family of operators defined in \eqref{eq:OperatoreBucoContinuo} with $M_0=\Id$ and $\mu$ be the unique invariant measure of the system \eqref{eq:system1} with density $f_0\in\Bexp(\R^d)$ such that \eqref{eq:VelocitaConvergenzaPalleEps} holds. Then, for every $f\in \Bexp(\R^d)$, we have
    \[
        \lim_{\xi\to 0^+}\left\|\frac{(M_\xi-M_0)f}{\xi}-\left(-\frac{\tau\, f(x_0)}{f_0(x_0)}\, \delta_{x_0}\right)\right\|_{KR}=0\ .
    \]

     \begin{myproof}
        For $\xi>0$, we have:
        \[
            \frac{(M_\xi-M_0)(f)(x)}{\xi} = -\frac{1}{\xi}\, f(x) \mathds{1}_{B_{r(\xi)}}(x)\ ,
        \]
        which belongs to $L^1(\R^d)$. We identify it with the signed measure $\nu_\xi := -\frac{1}{\xi}\, f(x) \mathds{1}_{B_{r(\xi)}}(x)\, dx$ and let us consider $\nu := -\frac{\tau\, f(x_0)}{f_0(x_0)}\, \delta_{x_0}$. We must show that
        \[
            \left\|\frac{(M_\xi-M_0)f}{\xi}-\left(-\frac{\tau\, f(x_0)}{f_0(x_0)}\, \delta_{x_0}\right)\right\|_{KR}=\sup \left\{ \int_{\R^d} g\, d(\nu_\xi - \nu) : \|g\|_\infty \leq 1,\, \Lip(g) \leq 1 \right\} \to 0\ ,
        \]
        as $\xi \to 0^+$.
        
        Let us define 
        \[
            c := \frac{m(B_{r(\xi)})}{\xi}=\frac{\tau}{f_0(x_0)}\ .
        \]
        For any test function $g:\R^d\to\R$ such that $\|g\|_\infty \leq 1$ and $\Lip(g) \leq 1$, we get:
        \begin{align*}
            \left| \int_{\R^d} g(x)\, d(\nu_\xi - \nu)(x) \right| &=\left|-\frac{1}{\xi}\int_{B_{r(\xi)}}f(x)g(x)\, dx+\frac{\tau\, f(x_0)g(x_0)}{f_0(x_0)}\right|= \left| \frac{c}{m(B_{r(\xi)})} \int_{B_{r(\xi)}} f(x)g(x)\, dx - c\, f(x_0)g(x_0) \right|\\
            &\leq\frac{c}{m(B_{r(\xi)})} \int_{B_{r(\xi)}} |f(x)g(x) - f(x_0)g(x_0)|\, dx\ .
        \end{align*}
        By adding and subtracting $f(x_0)g(x)$ and by using $\|g\|_\infty \leq 1$ and $\Lip(g) \leq 1$ (which implies $|g(x) - g(x_0)| \leq |x - x_0| \leq r(\xi)$ for $x \in B_{r(\xi)}$), we obtain:
        \[
            |f(x)g(x) - f(x_0)g(x_0)|\leq |f(x) - f(x_0)| |g(x)| + |f(x_0)| |g(x) - g(x_0)|\leq |f(x) - f(x_0)| + |f(x_0)|\, r(\xi)\ .
        \]
        Averaging this bound over the ball yields:
        \[
            \frac{c}{m(B_{r(\xi)})} \int_{B_{r(\xi)}} |f(x)g(x) - f(x_0)g(x_0)|\, dx \leq \frac{c}{m(B_{r(\xi)})} \int_{B_{r(\xi)}} |f(x) - f(x_0)|\, dx + c\, |f(x_0)|\, r(\xi)\ .
        \]
        Because $f \in \Bexp(\R^d) \subset C^1(\R^d)$ is continuous at $x_0$, the integral average in the right-hand side above goes to zero and $r(\xi) \to 0$ as $\xi \to 0^+$. Thus, the entire bound vanishes as $\xi \to 0^+$ uniformly for all valid $g$, concluding the proof.
    \end{myproof}
\end{prop}

Next, we compute the general derivative operator in a small neighborhood of zero of $M_\eps$ in the Kantorovich--Rubinstein norm. 

\begin{prop}\label{prop: ConvergenzaWassBucoFuoriZero}
    Let $x_0\in\R^d$ and $(M_\eps)_{\eps\in I_0^+}$ be the family of operators defined in \eqref{eq:OperatoreBucoContinuo}. Then, for every $f\in \Bexp(\R^d)$ and every $\eps\in I_0^+$, we have
    \[
        \lim_{\xi\to 0}\left\|\frac{(M_{\eps+\xi}-M_\eps)f}{\xi}-\left(-r'(\eps)f\,  \mathcal{H}^{d-1}_{|\de B_{r(\eps)}}\right)\right\|_{KR}=0\ ,
    \]
    where $\mathcal{H}^{d-1}_{|\de B_{r(\eps)}}$ denotes the $(d-1)$-dimensional Hausdorff measure restricted to the sphere centered at $x_0$ of radius $r(\eps)$.
    
    \begin{myproof}
        We consider the right limit as $\xi\to 0^+$ (the left limit as $\xi \to 0^-$ follows identically via the substitution $\delta = -\xi > 0$, replacing the annulus involved below with $B_{r(\eps)} \smallsetminus B_{r(\eps-\delta)}$). We note that
        \[
            \frac{(M_{\eps+\xi}-M_\eps)(f)(x)}{\xi}=-\frac{1}{\xi}\, f(x) \mathds{1}_{B_{r(\eps+\xi)}\smallsetminus B_{r(\eps)}}(x)\ ,
        \]
        which belongs to $L^1(\R^d)$, thus can be identified with the signed measure $\nu_{\eps,\xi} := -\frac{1}{\xi}\, f(x) \mathds{1}_{B_{r(\eps+\xi)}\smallsetminus B_{r(\eps)}}(x)\, dx$. We must show its Kantorovich--Rubinstein distance to $\nu_\eps := -r'(\eps) f(x)\, d\mathcal{H}^{d-1}_{|\de B_{r(\eps)}}(x)$ vanishes.
    
        For any test function $g:\R^d\to\R$ with $\|g\|_\infty\leq 1$ and $\Lip(g)\leq 1$, we apply the coarea formula to the integral against $\nu_{\eps,\xi}$ to obtain:
         \[
            \int_{\R^d} g(x)\, d\nu_{\eps,\xi}(x) = -\frac{1}{\xi}\int_{B_{r(\eps+\xi)}\smallsetminus B_{r(\eps)}}f(x)g(x)\, dx= -\frac{1}{\xi} \int_{r(\eps)}^{r(\eps+\xi)} G_g(s)\, ds\ ,
        \]
        where we set, for each $s>0$,
        \[
            G_g(s) := \int_{\de B_s} f(x)g(x)\, d\mathcal{H}^{d-1}(x)\ .
        \]
        In particular, the integral against the target measure $\nu_\eps$ is:
        \[
           \int_{\R^d}g(x)\, d\nu_\eps(x)=-r'(\eps)\int_{\de B_{r(\eps)}} f(x)g(x)\, d\mathcal{H}^{d-1}(x)=-r'(\eps)\, G_g(r(\eps))\ .
        \]
        If we add and subtract $\frac{r(\eps+\xi)-r(\eps)}{\xi} G_g(r(\eps))$, then $\left| \int_{\R^d} g(x)\, d(\nu_{\eps,\xi} - \nu_\eps)(x) \right|$ can be bounded by:
        \[
            \left| \frac{r(\eps+\xi)-r(\eps)}{\xi} - r'(\eps) \right| |G_g(r(\eps))| + \frac{r(\eps+\xi)-r(\eps)}{\xi} \left(\frac{1}{r(\eps+\xi)-r(\eps)} \int_{r(\eps)}^{r(\eps+\xi)} |G_g(s) - G_g(r(\eps))|\, ds\right)\ .
        \]
        Since $\|g\|_\infty \leq 1$, the first term is bounded by 
        \[
            \left|\frac{r(\eps+\xi)-r(\eps)}{\xi} - r'(\eps)\right| \|f\|_\infty\, \mathcal{H}^{d-1}(\de B_{r(\eps)})\ ,
        \]
        which vanishes as $\xi \to 0^+$, independently of $g$. For the second term, we show $G_g(s) \to G_g(r(\eps))$ uniformly over the class of test functions $g$. Indeed, using the radial projection $T_s(y) := x_0 + \frac{s}{r(\eps)}(y-x_0)$, we pull the integral back to the fixed sphere $\de B_{r(\eps)}$:
        \[
            G_g(s) = \int_{\de B_{r(\eps)}} f(T_s(y))g(T_s(y)) \left(\frac{s}{r(\eps)}\right)^{d-1} d\mathcal{H}^{d-1}(y)\ .
        \]
        Because $f\in \Bexp(\R^d)\subset C^1(\R^d)$, it is uniformly continuous on compact neighborhoods of $\de B_{r(\eps)}$. Moreover, the family of test functions $g$ is uniformly bounded and uniformly equicontinuous ($\Lip(g) \leq 1$). Since $T_s(y) \to y$ uniformly as $s \to r(\eps)$, the composite integrand converges to $f(y)g(y)$ uniformly in $y$ and uniformly over all admissible $g$. 
        Thus, $\sup_g |G_g(s) - G_g(r(\eps))| \to 0$ as $s \to r(\eps)$, causing the integral average over the shrinking interval $[r(\eps), r(\eps+\xi)]$ to vanish uniformly as $\xi \to 0^+$ (while the factor before the parenthesis converges to $r'(\eps)$). This concludes the proof.
    \end{myproof}
\end{prop}

\subsubsection*{Uniform continuity of the operator family}

The following integral will be useful for the computations below.

\begin{lem}\label{lem: integralacci Deboli}
    Let $\alpha,c>0$. Then, for each $x\in\R^d$, the integral
    \[
        I_{\alpha,c}(x):=\int_{\R^d}\rho_{\alpha,\exp}(y)\, e^{-c|x-y|^2}\, dy
    \]
    is finite.
    
    \begin{myproof}
        We first remark that $\alpha\, (1+|y|^2)^{\frac{1}{2}} \leq \alpha\, (1+|y|)$. Using the inequality $-|x-y|^2 \leq -|x|^2 + 2\, |x||y| - |y|^2$, we can bound the exponent in $I_{\alpha,c}(x)$ as follows:
        \[
            \alpha\, (1+|y|^2)^{\frac{1}{2}} - c\, |x-y|^2 \leq -c\, |y|^2 + (\alpha+2c\, |x|)|y| - c\, |x|^2 + \alpha \leq -\frac{c}{2}\, |y|^2 + \frac{1}{2c}\, (\alpha+2c|x|)^2+\alpha-c\, |x|^2\ , 
        \]
        where we applied Young's inequality to the last passage. Let $K(x) := \frac{1}{2c}\, (\alpha+2c|x|)^2 - c\, |x|^2 + \alpha$. Then
        \[
            I_{\alpha,c}(x) \leq \int_{\R^d} e^{-\frac{c}{2}|y|^2 + K(x)}\, dy = e^{K(x)} \left(\frac{2\pi}{c}\right)^{\frac{d}{2}}\ ,
        \]
        which is finite for every $x \in \R^d$.
    \end{myproof}
\end{lem}

We begin by establishing the continuity of the parameterization. By exploiting the analytic bounds of the transition density over the shrinking annuli between slightly perturbed balls, we show that $\tilde{\cl}_{h,\eps}$ depends continuously on $\eps$ in the operator norm.

\begin{prop}\label{prop: FamigliaContinua Weak}
    The following map:
    \[
        \funcdef{I_0}{\Ll(\Bexp(\R^d))}
        {\eps}{\tilde{\cl}_{h,\eps}}
    \]
    is continuous, i.e., for every $\eps_0\in I_0$,
    \[
        \lim_{\eps\to\eps_0}\|\tilde{\cl}_{h,\eps}-\tilde{\cl}_{h,\eps_0}\|_{\Ll(\Bexp)}=0\ .
    \]

    \begin{myproof}
        Let us define the annulus $A_{\eps,\eps_0} := B_{r(|\eps|)} \triangle B_{r(|\eps_0|)}$, where $\triangle$ denotes the symmetric difference of the two balls. By \eqref{eq:DefTransOpContinuo}, for $\eps$ sufficiently close to $\eps_0$ (sharing the same sign if $\eps_0 \neq 0$), the absolute pointwise difference between the operators reduces to an integral over $A_{\eps,\eps_0}$. We note that if $\eps_0 = 0$, this integration domain is simply the ball $B_{r(|\eps|)}$ rather than an annulus. Therefore, in all cases, we can uniformly bound the difference and its spatial gradient as follows:
        \[
            |(\tilde{\cl}_{h,\eps} - \tilde{\cl}_{h,\eps_0})(f)(y)| \leq \int_{A_{\eps,\eps_0}} S_h(x,y)|f(x)|\, dx \leq \|f\|_{\infty} \int_{A_{\eps,\eps_0}} S_h(x,y)\, dx\ ,
        \]
        \[
            |\nabla_y (\tilde{\cl}_{h,\eps} - \tilde{\cl}_{h,\eps_0})(f)(y)| \leq \int_{A_{\eps,\eps_0}} |\nabla_y S_h(x,y)f(x)|\, dx \leq \|f\|_{\infty} \int_{A_{\eps,\eps_0}} |\nabla_y S_h(x,y)|\, dx\ .
        \]
        Let $\bar{r} > 0$ be a fixed radius such that $A_{\eps,\eps_0} \subset \overline{B_{\bar{r}}(x_0)}$ for all $\eps$ near $\eps_0$. Applying the Gaussian bounds from Theorem \ref{teo: MPZ Weak} with $T=h$, for $x \in \overline{B_{\bar{r}}(x_0)}$ and using $|x-y|^2\geq \frac{1}{2}\, |x_0-y|^2-|x_0-x|^2$, we find uniform bounds over the domain of integration:
        \[
            \sup_{x \in A_{\eps,\eps_0}} S_h(x,y) \leq C_0\, e^{\lambda_0\|b\|_{\infty}^2T}\, h^{-\frac{d}{2}}\, e^{\frac{\lambda_0}{2h}\bar{r}^2}\, e^{-\frac{\lambda_0}{4h}|x_0-y|^2} =: \Lambda_0\, e^{-\frac{\lambda_0}{4h}|x_0-y|^2} 
        \]
        \[
            \sup_{x \in A_{\eps,\eps_0}} |\nabla_y\, S_h(x,y)| \leq C_1\, e^{\lambda_1\|b\|_\infty^2 T}\, h^{-\frac{1}{2}(d+1)}\, e^{\frac{\lambda_1}{2h}\bar{r}^2}\, e^{-\frac{\lambda_1}{4h}|x_0-y|^2} =: \Lambda_1\, e^{-\frac{\lambda_1}{4h}|x_0-y|^2}\ .
        \]
        Bounding the integrals above by these suprema times the Lebesgue measure $m(A_{\eps,\eps_0})$, we immediately obtain pointwise estimates:
        \[
            |(\tilde{\cl}_{h,\eps} - \tilde{\cl}_{h,\eps_0})(f)(y)| \leq \|f\|_\infty\, \Lambda_0\, e^{-\frac{\lambda_0}{4h}|x_0-y|^2}\, m(A_{\eps,\eps_0})\ ,
        \]
        \[
            |\nabla_y (\tilde{\cl}_{h,\eps} - \tilde{\cl}_{h,\eps_0})f(y)|\leq \|f\|_\infty\, \Lambda_1\, e^{-\frac{\lambda_1}{4h}|x_0-y|^2}\, m(A_{\eps,\eps_0})\ .
        \]
        We now compute the $\Bexp$-norm components. By bounding the exponential tails by $1$, we get the $L^\infty$-norms:
        \[
            \|(\tilde{\cl}_{h,\eps} - \tilde{\cl}_{h,\eps_0})f\|_\infty \leq \Lambda_0\, \|f\|_\infty\, m(A_{\eps,\eps_0}) , \qquad
            \|\nabla (\tilde{\cl}_{h,\eps} - \tilde{\cl}_{h,\eps_0})f\|_\infty \leq \Lambda_1\, \|f\|_\infty\, m(A_{\eps,\eps_0})\ ,
        \]
        thus
        \[
            \|(\tilde{\cl}_{h,\eps}-\tilde{\cl}_{h,\eps_0})f\|_{C^1}\leq \|f\|_\infty\, (\Lambda_0+\Lambda_1)\, m(A_{\eps,\eps_0})\ .
        \]
        Multiplying $|(\tilde{\cl}_{h,\eps}-\tilde{\cl}_{h,\eps_0})(f)(y)|$ by $\rho_{\alpha,\exp}(y)$ and integrating over $\R^d$ yields the $L^1_{\alpha,\exp}$-norm and via the previous estimate and Lemma \ref{lem: integralacci Deboli} we have:
        \[
            \|(\tilde{\cl}_{h,\eps}-\tilde{\cl}_{h,\eps_0})f\|_{L^1_{\alpha,\exp}} \leq\|f\|_{\infty}\, \Lambda_0\left(\int_{\R^d}\rho_{\alpha,\exp}(y)\, e^{-\frac{\lambda_0}{4h}|x_0-y|^2}\, dy\right) m(A_{\eps,\eps_0}) \leq \|f\|_\infty\, \Lambda_0\,  I_{\alpha,\frac{\lambda_0}{4h}}(x_0)\, m(A_{\eps,\eps_0})\ .
        \]
        Summing these components and using $\|f\|_\infty \leq \|f\|_{\Bexp}$, we conclude:
        \[
             \|(\tilde{\cl}_{h,\eps}-\tilde{\cl}_{h,\eps_0})f\|_{\Bexp} \leq \left(\Lambda_0\left(1+I_{\alpha,\frac{\lambda_0}{4h}}(x_0)\right)+\Lambda_1\right) \|f\|_{\Bexp}\, m(A_{\eps,\eps_0})\ .
        \]
        Let us set the constant $\bar{K} := \Lambda_0 \left(1 + I_{\alpha,\frac{\lambda_0}{4h}}(x_0)\right) + \Lambda_1$. Taking the supremum over $\|f\|_{\Bexp} \leq 1$ yields:
        \[
            \|\tilde{\cl}_{h,\eps}-\tilde{\cl}_{h,\eps_0}\|_{\Ll(\Bexp)} \leq \bar{K}\, m(A_{\eps,\eps_0})\ .
        \]
        Since $r(\eps)$ is a continuous function, $m(A_{\eps,\eps_0}) = |m(B_{r(|\eps|)}) - m(B_{r(|\eps_0|)})| \to 0$ as $\eps \to \eps_0$, which completes the proof.
    \end{myproof}
\end{prop}

\subsubsection*{Pointwise differentiability in the \texorpdfstring{$\Bexp$}{B-exp}-topology}\label{subsec: Derivate Transop}

Here we determine the pointwise derivative operator in the $\Bexp$-topology, i.e., we find an operator $\tilde{\cl}'_h(\eps)$ such that 
\[
    \lim_{\xi\to 0}\left\|\frac{(\tilde{\cl}_{h,\eps+\xi}-\tilde{\cl}_{h,\eps})f}{\xi}-\tilde{\cl}'_h(\eps)f\right\|_{\Bexp}=0\ .
\]
The central strategy here is to apply the smoothing properties of the unperturbed transfer operator $\tilde{\cl}_{h,0}$ to upgrade the Kantorovich--Rubinstein convergence established in Subsection \ref{subsec: Derivata operatore intermedio} into convergence in the strong $\Bexp$-norm.

\smallskip With a slight abuse of notation, we still denote by $\tilde{\cl}_{h,0}$ the transfer operator acting on Borel signed measures with bounded variation; namely, for each $\eta\in SM(\R^d)$ (thus with $\|\eta\|=|\eta|(\R^d)<+\infty$), we set
\begin{equation}\label{eq:TransNotPerturbedOnMeasures}
    (\tilde{\cl}_{h,0}\, \eta)(y):=\int_{\R^d}S_h(x,y)\, d\eta (x)\ .
\end{equation}
As already pointed out before Definition \ref{def: TransOpGeneralizzato}, in this way we naturally extend Definition \ref{def: TransOp}.

\smallskip Before computing the derivative operator, we need the following result, which ensures regularization properties of $\tilde{\cl}_{h,0}$ from the Kantorovich--Rubinstein norm to the $\Bexp$-norm, passing through the $L^1_{\alpha,\exp}$-norm.

\begin{prop}\label{prop: StimaTransOp Forte-Wasserstein}
    Let $h>0$, $x_0\in\R^d$ and $\mu,\nu\in SM(\R^d)$ such that their supports are contained in a closed ball $B$ centered at $x_0$ with radius $r>0$ and let $\tilde{\cl}_{h,0}$ be the transfer operator defined in \eqref{eq:TransNotPerturbedOnMeasures}. Then there exists a constant $\tilde{C}>0$ (depending on $h,d,\lambda_s,K,r,\alpha,\|b\|_{\infty}$ and $x_0$) such that
    \[
        \|\tilde{\cl}_{h,0}(\mu-\nu)\|_{\Bexp}\leq \tilde{C}\, \|\mu-\nu\|_{KR}\ .
    \]

    \begin{myproof}
        We first establish pointwise bounds for $\tilde{\cl}_{h,0}(\mu-\nu)(y)$ and its gradient. Let us consider the functions $k_y(x) := S_h(x,y)$ and $q_y(x) := \nabla_y\, S_h(x,y)$. Since $\mu-\nu$ has finite total variation and by Theorem \ref{teo: MPZ Weak} the map $y\mapsto \nabla_y\, S_h(x,y)$ is continuous and uniformly bounded on $B\times\R^d$, the dominated convergence theorem justifies differentiating under the integral sign and shows that $y\mapsto \nabla_y\, \tilde{\cl}_{h,0}(\mu-\nu)(y)$ is continuous on $\R^d$, thus $\tilde{\cl}_{h,0}(\mu-\nu)\in C^1_b(\R^d)$. By the definition of the Kantorovich--Rubinstein distance \eqref{eq: def Kantorovich--Rubinstein distance}, since the measures are supported in $B$, we have for any $y \in \R^d$:
        \begin{equation}\label{eq: KR_bound_func}
            |\tilde{\cl}_{h,0}(\mu-\nu)(y)| \leq \max\left(\|k_y\|_{\infty,B},\, \Lip(k_{y|B})\right) \|\mu-\nu\|_{KR}\ ,
        \end{equation}
        and
        \begin{equation}\label{eq: KR_bound_grad}
            |(\nabla_y \tilde{\cl}_{h,0}(\mu-\nu)(y))_i|\leq \max\left(\|q_y^i\|_{\infty,B},\, \Lip(q_{y|B}^i)\right) \|\mu-\nu\|_{KR} \ .
        \end{equation}
        where $q_y^i$ is the $i$-th component of $q_y$ (we are using the fact that any bounded Lipschitz function on $B$ extends to $\R^d$ with the exact same $\sup$-norm and Lipschitz constant). We highlight that in this proof we need the mixed derivative estimate from Appendix \ref{app: StimaD derivata mista} because by the convexity of $B$ and the mean value theorem, we can bound $\Lip(q_{y|B}^i)\leq \|\nabla_x\nabla_y\, S_h(\, \cdot\, , y)\|_{\infty,B}$.
        
        By Theorem \ref{teo: MPZ Weak} with  $T\geq h$, and the inequality $|x-y|^2\geq \frac{1}{2}\, |x_0-y|^2-|x_0-x|^2$, we evaluate these suprema over $x \in B$. For the function $k_y$, we have:
        \[
            \|k_y\|_{\infty,B} \leq C_0\, e^{\lambda_0\|b\|_{\infty}^2T}\, h^{-\frac{d}{2}}\, e^{\frac{\lambda_0}{2h}r^2}\, e^{-\frac{\lambda_0}{4h}|x_0-y|^2} \ ,\qquad
            \Lip(k_{y|B}) \leq C_1\, e^{\lambda_1\|b\|_\infty^2 T}\, h^{-\frac{1}{2}(d+1)}\, e^{\frac{\lambda_1}{2h}r^2}\, e^{-\frac{\lambda_1}{4h}|x_0-y|^2}\ .
        \]
        Let us define $\Lambda_s := C_s\, e^{\lambda_s\|b\|_\infty^2 T}\, h^{-\frac{1}{2}(d+s)}\, e^{\frac{\lambda_s}{2h}r^2}$ for $s \in \{0,1\}$. If we set $\Lambda_\star:=\max\{\Lambda_0,\Lambda_1\}$ and $\lambda_\star:=\min\{\lambda_0,\lambda_1\}$, then for each $s\in\{0,1\}$ and $y\in\R^d$,
        \begin{equation}\label{eq: Pointwise_Func}
            |\tilde{\cl}_{h,0}(\mu-\nu)(y)| \leq \Lambda_\star\, e^{-\frac{\lambda_\star}{4h}|x_0-y|^2} \|\mu-\nu\|_{KR}\ .
        \end{equation}
        
        Similarly, for the gradient components $q_y^i$, we apply Theorem \ref{teo: MPZ Weak} to $\nabla_y\, S_h$ and $\nabla_x\nabla_y\ S_h$. Bounding the exponential tail by $1$ for simplicity, we obtain:
        \[
            \|q_y^i\|_{\infty,B} \leq C_1\, e^{\lambda_1\|b\|_\infty^2\, T} h^{-\frac{1}{2}(d+1)} \ ,\qquad 
            \Lip(q_{y|B}^i) \leq C_2\, e^{\lambda_2\|b\|_\infty^2\, T} h^{-\frac{1}{2}(d+2)}
        \]
        Let $v \in \{1,2\}$ be the index that maximizes this bound in \eqref{eq: KR_bound_grad} and define $\tilde{\Lambda}_v := C_v\, e^{\lambda_v|b|_\infty^2 T}\, h^{-\frac{1}{2}(d+v)}$.
        Taking the Euclidean norm over the $d$ components yields:
        \begin{equation}\label{eq: Pointwise_Grad}
            |\nabla_y\, \tilde{\cl}_{h,0}(\mu-\nu)(y)| \leq d^{\frac{1}{2}} \, \tilde{\Lambda}_v\, 
            \|\mu-\nu\|_{KR}\ .
        \end{equation}

        We now compute the three required norms. By bounding the exponential term in \eqref{eq: Pointwise_Func} by $1$, we immediately obtain the $L^\infty$-norm:
        \[
            \|\tilde{\cl}_{h,0}(\mu-\nu)\|_{\infty} \leq \Lambda_\star \|\mu-\nu\|_{KR}\ ,\qquad \|\nabla\, \tilde{\cl}_{h,0}(\mu-\nu)\|_\infty\leq d^{\frac{1}{2}}\, \tilde{\Lambda}_v\, \|\mu-\nu\|_{KR}\ .
        \]
        Multiplying \eqref{eq: Pointwise_Func} by $\rho_{\alpha,\exp}(y)$ and integrating over $\R^d$ yields the $L^1_{\alpha,\exp}$-norm and via Lemma \ref{lem: integralacci Deboli} we get:
        \[
            \|\tilde{\cl}_{h,0}(\mu-\nu)\|_{L^1_{\alpha,\exp}} \leq \Lambda_\star\,  I_{\alpha,\frac{\lambda_\star}{4h}}(x_0)\, \|\mu-\nu\|_{KR}\ .
        \]
        
        Summing these three bounds gives the desired result, with the overall constant given by 
        \begin{equation}\label{eq: Costante Finale}
            \tilde{C} = \Lambda_\star\left(1 + I_{\alpha,\frac{\lambda_\star}{4h}}(x_0)\right) + d^{\frac{1}{2}}\, \tilde{\Lambda}_v\ .
        \end{equation}
        The proof is done.
    \end{myproof}
\end{prop}

\begin{prop}\label{prop: DerivataTransOp Weak}
    Let $x_0\in\R^d$ and $(\tilde{\cl}_{h,\eps})_{\eps\in I_0}$ be the family of transfer operators defined in \eqref{eq:DefTransOpContinuo}. With a slight abuse of notation, let $\tilde{\cl}_{h,0}$ still be the transfer operator acting on measures, defined in \eqref{eq:TransNotPerturbedOnMeasures} and $\mu$ be the unique invariant measure of the system \eqref{eq:system1} with density $f_0\in\Bexp(\R^d)$ such that \eqref{eq:VelocitaConvergenzaPalleEps} holds. Then, for each $f\in \Bexp(\R^d)$, we have
    \begin{equation}\label{eq: Limite1}
        \lim_{\xi\to 0}\left\|\frac{(\tilde{\cl}_{h,\xi}-\tilde{\cl}_{h,0})f}{\xi}-\tilde{\cl}_{h,0}\left(-\frac{\tau\, f(x_0)}{f_0(x_0)}\, \delta_{x_0}\right)\right\|_{\Bexp}=0\ .
    \end{equation} 
    Moreover, for each $f\in \Bexp(\R^d)$ and  $\eps\in I_0\smallsetminus\{0\}$, we have
    \begin{equation}\label{eq: Limite2}
        \lim_{\xi\to 0}\left\|\frac{(\tilde{\cl}_{h,\eps+\xi}-\tilde{\cl}_{h,\eps})f}{\xi}-\tilde{\cl}_{h,0}\left(-\sgn(\eps)\, r'(\eps)f\, \mathcal{H}^{d-1}_{|\de B_{r(\eps)}}\right)\right\|_{\Bexp}=0\ .
    \end{equation}
    
    \begin{myproof}
        $(A)$ \textit{Derivative at $\eps=0$}. We begin studying the case $\xi>0$. Recalling \eqref{eq: DefinizionTransFinaleComeComposizione}, we have
        \[
            \frac{(\tilde{\cl}_{h,\xi}-\tilde{\cl}_{h,0})f}{\xi}=\frac{\tilde{\cl}_{h,0}(M_\xi f-M_0f)}{\xi}=\tilde{\cl}_{h,0}\left(\frac{(M_\xi-M_0)f}{\xi}\right)\ .
        \] 
        Let $\nu_0 := -\frac{\tau\, f(x_0)}{f_0(x_0)}\, \delta_{x_0}$ and $\eta_\xi := \frac{(M_\xi - M_0)f}{\xi}$. By Propositions \ref{prop:ConvergenzaWassBuco} and \ref{prop: StimaTransOp Forte-Wasserstein}, we conclude that
        \[
            \|\tilde{\cl}_{h,0}(\eta_\xi-\nu_0)\|_{\Bexp}\leq \tilde{C}\, \|\eta_\xi-\nu_0\|_{KR}\to 0
        \]
        as $\xi\to 0^+$. Now we study the case $\xi<0$. By \eqref{eq: DefinizionTransFinaleComeComposizione}, setting $\delta=-\xi>0$, we have
        \[
            \frac{(\tilde{\cl}_{h,\xi}-\tilde{\cl}_{h,0})f}{\xi}=\frac{\tilde{\cl}_{h,0}(2M_0-M_{-\xi}-M_0)f}{\xi}=\frac{\tilde{\cl}_{h,0}(M_0 -M_{-\xi})f}{-(-\xi)}=\tilde{\cl}_{h,0}\left(\frac{(M_\delta-M_0)f}{\delta}\right)\ 
        \]
        and the conclusion is analogous by taking the limit as $\delta\to 0^+$. The explicit pointwise limit is:
        \[
            (\tilde{\cl}_{h,0}\, \nu_0)(y)=\int_{\R^d} S_h(x,y)\, d\nu_0(x)=-\frac{\tau\, f(x_0)S_h(x_0,y)}{f_0(x_0)} \ .
        \]
        This establishes \eqref{eq: Limite1}.

        \smallskip $(B)$ \textit{Derivative at $\eps \neq 0$}. We begin studying the case $\eps\in I_0^+$. By \eqref{eq: DefinizionTransFinaleComeComposizione},
        \[
            \frac{(\tilde{\cl}_{h,\eps+\xi}-\tilde{\cl}_{h,\eps})f}{\xi}=\frac{\tilde{\cl}_{h,0}(M_{\eps+\xi} f-M_\eps\, f)}{\xi}=\tilde{\cl}_{h,0}\left(\frac{(M_{\eps+\xi}-M_\eps)f}{\xi}\right)\ .
        \]
        Let $\nu_\eps^+ := - r'(\eps)f\, \mathcal{H}^{d-1}_{|\de B_{r(\eps)}}$ and $\eta_{\eps,\xi}:= \frac{(M_{\eps+\xi} - M_\eps)f}{\xi}\, dx$. By Propositions \ref{prop: ConvergenzaWassBucoFuoriZero} and \ref{prop: StimaTransOp Forte-Wasserstein}, we conclude that
        \[
            \|\tilde{\cl}_{h,0}(\eta_{\eps,\xi}-\nu_\eps^+)\|_{\Bexp}\leq\tilde{C}\, \|\eta_{\eps,\xi}-\nu_\eps^+\|_{KR}\to 0\ ,
        \]
        as $\xi\to 0$. Now we study the case $\eps\in I_0^-$. By \eqref{eq: DefinizionTransFinaleComeComposizione} and the symmetry $M_t = M_{-t}$, setting $\delta=-\xi$, we note that the difference quotient can be rewritten in terms of strictly positive indices:
        \[
            \frac{(\tilde{\cl}_{h,\eps+\xi}-\tilde{\cl}_{h,\eps})f}{\xi}=\frac{\tilde{\cl}_{h,0}(2M_0-M_{-\eps-\xi}-2M_0+M_{-\eps} )f}{\xi}=\tilde{\cl}_{h,0}\left(-\frac{(M_{-\eps-\xi}-M_{-\eps})f}{\xi}\right)=\tilde{\cl}_{h,0}\left(\frac{(M_{-\eps+\delta}-M_{-\eps})f}{\delta}\right)\ .
        \]
        Taking the limit as $\delta\to 0$ and by the previous case ($-\eps\in I_0^+$) yields the limit measure $\nu_\eps^-:=-r'(-\eps)f\, \mathcal{H}_{|\de B_{r(-\eps)}}^{d-1}=r'(\eps)f\, \mathcal{H}_{|\de B_{r(\eps)}}^{d-1}$, where we used that $r$ is an even function. 

        Both cases $\eps\in I_0^+$ and $\eps\in I_0^-$ can be unified by defining $\nu_\eps:=-\sgn(\eps)\, r'(\eps)f\, \mathcal{H}_{|\de B_{r(\eps)}}^{d-1}$, yielding the explicit pointwise limit: 
        \[
            (\tilde{\cl}_{h,0}\, \nu_\eps)(y)=\int_{\R^d}S_h(x,y)\, d\nu_\eps(x)=-\sgn(\eps)\, r'(\eps)\int_{\de B_{r(\eps)}}S_h(x,y)f(x)\, d\mathcal{H}^{d-1}(x) \ ,
        \]
        which establishes \eqref{eq: Limite2} and completes the proof.
    \end{myproof}
\end{prop}

We have just proved that the map $\eps\mapsto \tilde{\cl}_{h,\eps} f$ is differentiable for each $f\in \Bexp(\R^d)$ and $\eps\in I_0$. In particular,
for each $y\in\R^d$, its expression is given by:
\begin{equation}\label{eq:DerivataPuntuale TransOp}
    (\tilde{\cl}'_h(\eps)f)(y)=\begin{cases}
        -\dfrac{\tau\, f(x_0)S_h(x_0,y)}{f_0(x_0)} & \eps=0\\[.3cm]
        -\sgn(\eps)\, r'(\eps)\displaystyle\int_{\de B_{r(\eps)}}S_h(x,y)f(x)\, d\mathcal{H}^{d-1}(x) & \eps\in I_0\smallsetminus\{0\}
    \end{cases}\ .
\end{equation}
We need to check that we get bounded operators on $\Bexp(\R^d)$.

\begin{prop}\label{prop: DerivativeBounded}
    The derivative operator $\tilde{\cl}'_h(\eps)$ is bounded for each $\eps\in I_0$, i.e., 
    \[
            \|\tilde{\cl}'_h(\eps)\|_{\Ll(\Bexp)}=\sup_{\|f\|_{\Bexp}\leq 1}\|\tilde{\cl}'_h(\eps)f\|_{\Bexp}<+\infty\ .
    \]

    \begin{myproof}
        By Proposition \ref{prop: DerivataTransOp Weak}, we can express the action of the derivative operator as $\tilde{\cl}'_h(\eps) f=\tilde{\cl}_{h,0}\, \nu_{\eps,f}$, where $\tilde{\cl}_{h,0}$ is the transfer operator extended to measures defined in \eqref{eq:TransNotPerturbedOnMeasures} and $\nu_{\eps,f}$ is the corresponding finite signed measure supported on the compact ball $\overline{B_{r(\eps)}(x_0)}$ for every $\eps\in I_0$.
        In both cases (whether $\eps=0$ or $\eps\neq 0$), we can bound the total variation of this measure by a constant $K_0$ independent of $f$ and $\eps$. Specifically, for $\eps=0$, the measure is $\nu_{0,f}=-\frac{\tau\, f(x_0)}{f_0(x_0)}\, \delta_{x_0}$; for $\eps\neq 0$, the measure is $\nu_{\eps,f}=-\sgn(\eps)\, r'(\eps)f\, \mathcal{H}_{|\de B_{r(\eps)}}^{d-1}$. Because of the identity $|r'(\eps)|\, \mathcal{H}^{d-1}(\de B_{r(\eps)})=\frac{\tau}{f_0(x_0)}$, both cases yield: 
        \[
            \|\nu_{\eps,f}\|_{TV} \leq K_0\, \|f\|_\infty\leq K_0\, \|f\|_{\Bexp}\ ,
        \]
        where $K_0:=\frac{\tau}{f_0(x_0)}$.
        
        For every finite signed measure $\eta\in SM(\R^d)$, the definition of the Kantorovich--Rubinstein norm gives $\|\eta\|_{KR}\leq \|\eta\|_{TV}$. Because, for all $\eps\in I_0$, the measure $\nu_{\eps,f}$ is supported in the fixed closed ball $B_{\overline{r}}(x_0)$ where $\overline{r}:=\sup_{\eps\in I_0}|r(\eps)|<+\infty$, applying Proposition \ref{prop: StimaTransOp Forte-Wasserstein} with $\mu=\nu_{\eps,f}$ and $\nu=0$ gives
        \[
            \|\tilde{\cl}'_h(\eps)f\|_{\Bexp} = \|\tilde{\cl}_{h,0}\, \nu_{\eps,f}\|_{\Bexp} \leq \tilde{C}\|\nu_{\eps,f}\|_{KR}\leq \tilde{C}\, \|\nu_{\eps,f}\|_{TV} \leq \tilde{C} K_0\, \|f\|_{\Bexp}\ .
        \]
        Taking the supremum over all $f \in \Bexp(\R^d)$ with $\|f\|_{\Bexp} \leq 1$ yields $\|\tilde{\cl}'_h(\eps)\|_{\Ll(\Bexp)} \leq \tilde{C}K_0 < +\infty$, which completes the proof.
    \end{myproof}
\end{prop}

\subsubsection*{Continuous differentiability} 

Finally, to conclude the proof of $C^1$-regularity, we must verify that the map  $I_0\ni \eps \mapsto \tilde{\cl}'_h(\eps)\in\Ll(\Bexp(\R^d))$ is continuous. This successfully establishes that the family $(\tilde{\cl}_{h,\eps})_{\eps \in I_0}$ is of class $C^1(I_0,\Ll(\Bexp(\R^d)))$.

\begin{prop}\label{prop: Continuità Derivata Debole}
    The map 
    \[
        \funcdef{I_0}{\Ll(\Bexp(\R^d))}
        {\eps}{\tilde{\cl}'_h(\eps)}
    \]
    is continuous, i.e., for every $\eps_0\in I_0$,
    \[
        \lim_{\eps\to\eps_0}\|\tilde{\cl}'_h(\eps)-\tilde{\cl}'_h(\eps_0)\|_{\Ll(\Bexp)}=0\ .
    \]
    Therefore the family $(\tilde{\cl}_{h,\eps})_{\eps\in I_0}$ is of class $C^1(I_0, \Ll(\Bexp(\R^d)))$.

    \begin{myproof}
        By Proposition \ref{prop: DerivataTransOp Weak}, for any $\eps \in I_0$, we can express the derivative as $\tilde{\cl}'_h(\eps)f = \tilde{\cl}_{h,0}\, \nu_{\eps,f}$, where $\nu_{\eps,f}$ is the corresponding limit measure. Proposition \ref{prop: StimaTransOp Forte-Wasserstein} guarantees that for each $\eps$ near $\eps_0$, there is a uniform constant $\tilde{C}$ (independent of $f$, see \eqref{eq: Costante Finale}) such that:
        \[
            \|\tilde{\cl}'_h(\eps) - \tilde{\cl}'_h(\eps_0)\|_{\Ll(\Bexp)} = \sup_{\|f\|_{\Bexp}\leq 1} \|\tilde{\cl}_{h,0}(\nu_{\eps,f}-\nu_{\eps_0,f})\|_{\Bexp} \leq \tilde{C} \sup_{\|f\|_{\Bexp}\leq 1} \|\nu_{\eps,f} - \nu_{\eps_0,f}\|_{KR}\ .
        \]
        Then it suffices to show that the Kantorovich--Rubinstein distance vanishes uniformly for all $f$ in the unit ball of $\Bexp(\R^d)$.

        \smallskip\noindent $(A)$ \textit{Continuity at $\eps_0 \neq 0$}. 
        For $\eps$ sufficiently close to $\eps_0$ (sharing the same sign), both measures are surface integrals over spheres of radius $r(\eps)$ and $r(\eps_0)$,
        \[
            \int_{\R^d}g(x)\, d(\nu_{\eps,f}-\nu_{\eps_0,f})(x)=-\sgn(\eps)\, r'(\eps)\int_{\de B_{r(\eps)}}f(x)g(x)\, d\mathcal{H}^{d-1}(x)+\sgn(\eps_0)\, r'(\eps_0)\int_{\de B_{r(\eps_0)}}f(x)g(x)\, d\mathcal{H}^{d-1}(x)\ .
        \]
        Using the radial projection map $T(y) := x_0 + \frac{r(\eps)}{r(\eps_0)}(y-x_0)$, we can pull the integral for $\nu_{\eps,f}$ back to $\de B_{r(\eps_0)}$. By differentiating the explicit formula for $r(\eps)^d$ with respect to $\eps$ (see \eqref{eq:ChosenRadiusFunction}), we observe the identity $\sgn(\eps)\, r'(\eps)\, r(\eps)^{d-1} = \sgn(\eps_0)\, r'(\eps_0)\, r(\eps_0)^{d-1}$ (both are equal to $\frac{\tau}{f_0(x_0)\, \omega_{d-1}}$).
        This cancellation allows us to rewrite the difference over a single domain:
        \[
            \int_{\R^d}g(x)\, d(\nu_{\eps,f}-\nu_{\eps_0,f})(x) = -\sgn(\eps_0)\, r'(\eps_0)\int_{\de B_{r(\eps_0)}}\left(f(T(y))g(T(y))-f(y)g(y)\right)d\mathcal{H}^{d-1}(y)\ .
        \]
        By adding and subtracting $f(y)g(T(y))$ and applying the triangle inequality alongside the bounds $\|g\|_\infty, \Lip(g) \le 1$, we obtain:
        \[  
            \left|\int_{\R^d}g(x)\, d(\nu_{\eps,f}-\nu_{\eps_0,f})(x)\right| \leq |r'(\eps_0)| \int_{\de B_{r(\eps_0)}}\left(|f(T(y))-f(y)|+|r(\eps)-r(\eps_0)||f(y)|\right) d\mathcal{H}^{d-1}(y)\ .
        \]
        For every $f\in \Bexp(\R^d)$ such that $\|f\|_{\Bexp}\leq 1$ (thus $\|f\|_{C^1}\leq 1$), the integrand is bounded by $2\, |r(\eps)-r(\eps_0)|$. Taking the supremum over all test functions $g$ and then the supremum over all $f$ in the unit ball of $\Bexp(\R^d)$ yields:
        \[
            \sup_{\|f\|_{\Bexp}\leq 1}\|\nu_{\eps,f}-\nu_{\eps_0,f}\|_{KR}\leq |r'(\eps_0)|\left(2\, |r(\eps)-r(\eps_0)|\, \mathcal{H}^{d-1}(\de B_{r(\eps_0)})\right)=\frac{2\, \tau}{f_0(x_0)}\, |r(\eps)-r(\eps_0)|\ .
        \]
        The rightmost term vanishes as $\eps\to\eps_0$ because $r$ is continuous.

        \smallskip\noindent $(B)$ \textit{Continuity at $\eps_0 = 0$}. 
        For $\eps \to 0$, we compare the surface measure $\nu_{\eps,f} = -\sgn(\eps)\, r'(\eps)f\, \mathcal{H}^{d-1}_{|\de B_{r(\eps)}}$ to the Dirac measure $\nu_{0,f} = -\frac{\tau\, f(x_0)}{f_0(x_0)}\, \delta_{x_0}$. Differentiating the expression for $r(\eps)^d$ reveals that the unweighted mass of the surface measure is identically constant for all $\eps\neq 0$: $\sgn(\eps)\, r'(\eps)\, \mathcal{H}^{d-1}(\de B_{r(\eps)})=\frac{\tau}{f_0(x_0)}$. This identity allows us to write the Dirac evaluation as an integral over the boundary:
        \[
            \int_{\R^d}g(x)\, d(\nu_{\eps,f}-\nu_{0,f})(x)=-\sgn(\eps)\, r'(\eps)\int_{\de B_{r(\eps)}}f(x)g(x)\, d\mathcal{H}^{d-1}(x)+\frac{\tau\, f(x_0)g(x_0)}{f_0(x_0)}\ .
        \]
        Therefore 
        \[
            \left|\int_{\R^d}g(x)\, d(\nu_{\eps,f}-\nu_{0,f})(x)\right|\leq |r'(\eps)|\int_{\de B_{r(\eps)}}|f(x_0)g(x_0)-f(x)g(x)|\, d\mathcal{H}^{d-1}(x)\ .
        \]
        By adding and subtracting $f(x_0)g(x)$ and using the bounds for $f$ and $g$, we note that the integrand is bounded by $2\, r(\eps)$.  Then taking the supremum over all admissible test functions $g$ and all $f$ in the unit ball of $\Bexp(\R^d)$ yields:
        \[
            \sup_{\|f\|_{\Bexp}\leq 1}\|\nu_{\eps,f}-\nu_{0,f}\|_{KR} \leq |r'(\eps)|\, 2\, r(\eps)\, \mathcal{H}^{d-1}(\de B_{r(\eps)}) = \frac{2\, \tau}{f_0(x_0)}\, r(\eps)\ .
        \]
        The rightmost term converges to zero as $\eps\to 0$.

        \smallskip Lastly, we deduce that the map $\eps\mapsto\tilde{\cl}_{h,\eps}$ is differentiable in operator norm on $I_0$ and hence belongs to $C^1(I_0,\Ll(\Bexp(\R^d)))$. Indeed, fix $\eps\in I_0$ and $\xi\neq 0$ such that $\eps+\xi\in I_0$. Since $I_0\ni s\mapsto \tilde{\cl}_{h,s}f\in\Bexp(\R^d)$ is continuously differentiable for every $f\in \Bexp(\R^d)$ with derivative $\tilde{\cl}_h'(s)f$ by Proposition \ref{prop: DerivataTransOp Weak} and the continuity of $\tilde{\cl}'_h$ proved above, the fundamental theorem of calculus gives
        \[
            \frac{(\tilde{\cl}_{h,\eps+\xi}-\tilde{\cl}_{h,\eps})f}{\xi}-\tilde{\cl}_h'(\eps)f=\frac{1}{\xi}\int_\eps^{\eps+\xi}(\tilde{\cl}_h'(s)-\tilde{\cl}_h'(\eps))f\, ds\ .
        \]
        Taking the $\Bexp$-norm and then the supremum over all $f\in\Bexp(\R^d)$ with $\|f\|_{\Bexp}\leq 1$ yields
        \[
            \left\|\frac{\tilde{\cl}_{h,\eps+\xi}-\tilde{\cl}_{h,\eps}}{\xi}-\tilde{\cl}_h'(\eps)\right\|_{\Ll(\Bexp)}\leq\sup_{s\in\left[\min\{\eps, \eps+\xi\},\, \max\{\eps,\eps+\xi\}\right]}\|\tilde{\cl}_h'(s)-\tilde{\cl}_h'(\eps)\|_{\Ll(\Bexp)}\to 0\ ,
        \]
        as $\xi\to 0$. Therefore the map $I_0\ni\eps\mapsto\tilde{\cl}_{h,\eps}\in \Ll(\Bexp(\R^d))$ is continuously differentiable in $\Ll(\Bexp(\R^d))$, completing the proof.
    \end{myproof}
\end{prop}

\subsection{Verifying assumption \texorpdfstring{$(H2)$}{(H2)}}\label{subsec: Spectral properties unperturbed transop}

We show that the unperturbed transfer operator $\tilde{\cl}_{h,0}$ is quasi-compact when acting on $\Bexp(\R^d)$, which verifies hypothesis $(H2)$ of Theorem \ref{teo: TeoPerturbativo}.

\begin{prop}\label{prop:QuasiCompactnessForUnperturbedTransOpExp}
    The unperturbed transfer operator $\tilde{\cl}_{h,0}$ is quasi-compact on $\Bexp(\R^d)$.

    \begin{myproof}
        We apply Theorem \ref{teo: QC} to the Banach space $(\Bexp(\R^d),\|\cdot\|_{\Bexp})$, using the $L^1$-norm as the continuous semi-norm $\|\cdot\|'$. Since the exponential weight $\rho_{\alpha,\exp}(x) \geq 1$, we trivially have $\|f\|_{L^1}\leq \|f\|_{\Bexp}$ for every $f\in \Bexp(\R^d)$. 

        \smallskip \textit{Assumption (a)}. Let $B$ be the closed unit ball in $\Bexp(\R^d)$. By Proposition \ref{prop: CompactEmbedding SpazioForte Weak}, given $(f_n)_{n\geq 1}\subset B$, there exists a subsequence $(f_{n_k})_{k\geq 1}$ and $\tilde{f}\in L^1(\R^d)$ such that $\|f_{n_k}-\tilde{f}\|_{L^1}\to 0$ as $k\to+\infty$. Up to an a.e. convergent subsequence, by Fatou's lemma, $\|\tilde{f}\|_{L^1_{\alpha,\exp}}\leq \liminf_{k}\|f_{n_k}\|_{L^1_{\alpha,\exp}}\leq 1$, thus $\tilde{f}\in L^1_{\alpha,\exp}(\R^d)$. Let us define $g:=\tilde{\cl}_{h,0}\tilde{f}$. By Lemma \ref{lem: C1Estimate} and Lemma \ref{lem: LSY L12 Weak}, $g\in \Bexp(\R^d)$, and since $\tilde{\cl}_{h,0}$ is an $L^1$-contraction, $\|\tilde{\cl}_{h,0}f_{n_k}-g\|_{L^1}\leq \|f_{n_k}-\tilde{f}\|_{L^1}\to 0$
        as $k\to+\infty$, as wanted.
        
        \smallskip\textit{Assumption (b)}. As shown in Lemma \ref{lem:TransOpWeakContraction}, $\tilde{\cl}_{h,0}$ is a Markov operator and hence a contraction on $L^1(\R^d)$. Thus, the inequality $\|\tilde{\cl}_{h,0} f\|_{L^1}\leq M\|f\|_{L^1}$ holds with $M=1$. 
        
        \smallskip\textit{Assumption (c)}. Consider the non-zero bounded linear functional $\mu_0\in(\Bexp(\R^d))^\ast$ defined by $\langle\mu_0,f\rangle:=\int_{\R^d}f(x)\, dx$.  By Lemma \ref{lem:TransOpWeakContraction}, $\tilde{\cl}_{h,0}$ preserves the integral, thus $\langle(\tilde{\cl}_{h,0})^\ast\mu_0,f\rangle=\langle\mu_0,\tilde{\cl}_{h,0}f\rangle=\langle\mu_0,f\rangle$ for all $f\in\Bexp(\R^d)$. Thus $\lambda_0=1$ is in the spectrum of $\tilde{\cl}_{h,0}^\ast$, which implies that $\lambda_0\in\sigma(\tilde{\cl}_{h,0})$ and hence $r(\tilde{\cl}_{h,0}) \geq 1$. By the Lasota-Yorke inequality in Proposition \ref{prop: LYITransWeak}, we have $\|\tilde{\cl}_{h,0}^n\|_{\Ll(\Bexp)}\leq \tilde{A}\lambda^n+\tilde{B}$, which implies 
        \[
            r(\tilde{\cl}_{h,0}) \leq \lim_{n\to+\infty}(\tilde{A}\lambda^n+\tilde{B})^{\frac{1}{n}} \leq 1\ .
        \]
        Therefore, $r(\tilde{\cl}_{h,0})=1$.
        Furthermore, for any $k\geq 1$, the Lasota-Yorke inequality gives:
        \[
            \|\tilde{\cl}_{h,0}^k f\|_{\Bexp}\leq \tilde{A}\lambda^k\, \|f\|_{\Bexp} + \tilde{B}\|f\|_{L^1}\ .
        \]
        Since $\lambda \in (0,1)$, we can choose $k$ large enough so that $\tilde{A}\lambda^k < 1$. Setting $r = (\tilde{A}\lambda^k)^{\frac{1}{k}}$ and $R = \tilde{B}$, we satisfy $r < r(\tilde{\cl}_{h,0}) = 1$ and $\|\tilde{\cl}_{h,0}^kf\|_{\Bexp}\leq r^k\|f\|_{\Bexp}+R\|f\|_{L^1}$, completing the verification.
    \end{myproof}
\end{prop}

\begin{prop}\label{prop: AutovaloreUnperturbedTransOp}
    The unperturbed transfer operator $\tilde{\cl}_{h,0}$ has $\lambda_0=1$ as its simple and only eigenvalue on the unit circle, to which corresponds a unique invariant probability density $f_0\in\Bexp(\R^d)$.

    \begin{myproof}
        By Proposition \ref{prop:QuasiCompactnessForUnperturbedTransOpExp}, $\tilde{\cl}_{h,0}$ is quasi-compact on $\Bexp(\R^d)$ with spectral radius $r(\tilde{\cl}_{h,0})=1$ and $1\in\sigma(\tilde{\cl}_{h,0})$. Moreover, by definition, $\lambda_0=1$ is an eigenvalue of $\tilde{\cl}_{h,0}$. 

        \smallskip Now, let us show that $\lambda_0=1$ is the only eigenvalue on the unit circle. Suppose that $\lambda\in\C$ with $|\lambda|=1$ and $f\in\Bexp(\R^d)\smallsetminus\{0\}$ are such that $\tilde{\cl}_{h,0}f=\lambda f$. Then, for each $y\in\R^d$,
        \[
            |f(y)|=|\lambda f(y)|=|(\tilde{\cl}_{h,0}f)(y)|=\left|\int_{\R^d}S_h(x,y)f(x)\, dx\right|\leq \int_{\R^d}S_h(x,y)|f(x)|\, dx=(\tilde{\cl}_{h,0}|f|)(y)\ .
        \]
        Integrating both sides over $\R^d$ and using the integral preservation of $\tilde{\cl}_{h,0}$, given by Lemma \ref{lem:TransOpWeakContraction}, we obtain  
        \[
            0\leq \int_{\R^d}\left((\tilde{\cl}_{h,0}|f|)(y)-|f(y)|\right)dy=0\ .
        \]
        Since $\tilde{\cl}_{h,0}|f|-|f|\geq 0$ is continuous (by Lemma \ref{lem: C1Estimate}), and has zero integral, it vanishes identically. Thus $|(\tilde{\cl}_{h,0}f)(y)|=(\tilde{\cl}_{h,0}|f|)(y)$ for every $y\in\R^d$, i.e., 
        \[
            \left|\int_{\R^d} S_h(x,y)f(x)\, dx\right|=\int_{\R^d}S_h(x,y)|f(x)|\, dx\ .
        \]
        Fix any $y\in\R^d$. Because $S_h$ is strictly positive everywhere and $f$ is continuous and not identically zero, $|f|$ is positive on some open set, so $(\tilde{\cl}_{h,0}|f|)(y)>0$. Hence $(\tilde{\cl}_{h,0}f)(y)\neq 0$. Write the integral as follows: $(\tilde{\cl}_{h,0}f)(y):=I=e^{i\theta}|I|$ for some $\theta:=\arg(I)\in\R$. Then
        \[
            \int_{\R^d}S_h(x,y)\left(|f(x)|-\Re (e^{-i\theta} f(x))\right)dx=(\tilde{\cl}_{h,0}|f|)(y)-|(\tilde{\cl}_{h,0}f)(y)|=0\ .
        \]
        By continuity of $f$ and strict positivity of $S_h$, this implies that the integrand must vanish identically, i.e.,  $\Re(e^{-i\theta}f(x))=|f(x)|$, and hence $e^{-i\theta}f(x)=|f(x)|\geq 0$ for all $x\in\R^d$. Setting $g:=e^{-i\theta}f\geq 0$, we have $\lambda g(z)=(\tilde{\cl}_{h,0}g)(z)=\int_{\R^d}S_h(x,z)g(x)\, dx>0$ for each $z\in\R^d$. Since $g\geq 0$ is not identically zero, this forces $\lambda>0$, and combined with $|\lambda|=1$, we get $\lambda=1$. 

        \smallskip The argument above also shows that every non-zero eigenfunction for $\lambda_0=1$ is a non-zero complex multiple of a strictly positive function, and therefore has a non-zero integral. In particular, if $f_1,f_2\in\ker(\tilde{\cl}_{h,0}-\Id)$ are non-zero, normalizing them so that $\int_{\R^d} f_1(x)\, dx=\int_{\R^d} f_2(x)\, dx=1$, gives $f_1-f_2\in\ker(\tilde{\cl}_{h,0}-\Id)$ with $\int_{\R^d}(f_1(x)-f_2(x))\, dx=0$, which forces $f_1=f_2$. Hence $\ker(\tilde{\cl}_{h,0}-\Id)$ is one-dimensional, spanned by a unique strictly positive probability density $f_0\in\Bexp(\R^d)$.  

        \smallskip Lastly, to establish that the eigenvalue is simple, suppose $g\in \Bexp(\R^d)$ satisfies $(\tilde{\cl}_{h,0}-\Id)^2 g=0$ and set $v:=(\tilde{\cl}_{h,0}-Id)g$. Then $v\in\ker(\tilde{\cl}_{h,0}-\Id)=\spann(f_0)$. So $v(x)=c f_0(x)$ for some $c\in\C$. On the other hand, integral preservation implies 
        \[
            c=\int_{\R^d}c f_0(x)\, dx=\int_{\R^d}v(x)\, dx=\int_{\R^d}(\tilde{\cl}_{h,0}g-g)(x)\, dx=0\ ,
        \]
        whence $v=0$. Thus $\ker((\tilde{\cl}_{h,0}-\Id)^2)=\ker(\tilde{\cl}_{h,0}-\Id)$. Hence $\lambda_0=1$ is a simple dominating eigenvalue.
     \end{myproof}
\end{prop}

\begin{oss}\label{oss: MisuraInvarianteSDE}
    By Propositions \ref{prop: LYITransWeak} and \ref{prop: AutovaloreUnperturbedTransOp}, we note that $\cl_t f_0\in\Bexp(\R^d)$. Moreover, $\cl_t f_0$ is a probability density. Thanks to the semigroup property and the fact that $f_0$ is fixed by $\cl_h$, we have $\cl_h(\cl_t f_0)=\cl_t(\cl_h f_0)=\cl_t f_0$. Since $\cl_h$ has a unique invariant probability density in $\Bexp(\R^d)$, it must be $\cl_t f_0=f_0$ for each $t>0$. Therefore, for each $\phi\in L^\infty(\R^d)$, setting $\nu=f_0(x)\, dx$, by the duality relation \eqref{eq:duality}, we have
    \[
        \int_{\R^d}(P_t\phi)(x) d\nu(x)=\int_{\R^d}\phi(y)(\cl_t f_0)(y)\, dy=\int_{\R^d}\phi(y) f_0(y)\, dy=\int_{\R^d}\phi(y)\, d\nu(y)\ .
    \]
    Hence $\nu$ is a stationary measure for the system \eqref{eq:system1} and must coincide with $\mu$ by uniqueness. Indeed, if $\tilde{\mu}$ is any invariant probability measure of \eqref{eq:system1}, then $\tilde{\mu}(A)=\int_{\R^d}\left(\int_A S_h(x,y)\, dy\right) d\tilde{\mu}(x)$, so $\tilde{\mu}$ is absolutely continuous with density $g(y):=(\cl_h\tilde{\mu})(y)\in L^1(\R^d)$ satisfying $\cl_h g=g$. Setting $\tilde{f}=g-f_0\in L^1(\R^d)$, we have $\cl_h \tilde{f}=\tilde{f}$ and $\int_{\R^d} \tilde{f}(x)\, dx=0$. Since $\tilde{f}\in C^1_b(\R^d)$ by Lemma \ref{lem: C1Estimate}, the strict positivity argument of Proposition \ref{prop: AutovaloreUnperturbedTransOp} implies that $\tilde{f}$ has constant sign, thus $\tilde{f}=0$ a.e., and we can conclude that $\tilde{\mu}=\nu$.
\end{oss}

\subsection{Proof of Theorem \ref{teo:Teorema1}}\label{subsec: Final proof}

In this final section, we conclude with the proof of the EVL stated in Theorem \ref{teo:Teorema1}.

\begin{oss}\label{oss: Conseguenze Hennion Herve}
    Having established the $C^1$-regularity of the transfer operator family $(\tilde{\cl}_{h,\eps})_{\eps\in I_0}$ and the quasi-compactness of the unperturbed operator $\tilde{\cl}_{h,0}$, all hypotheses of Theorem \ref{teo: TeoPerturbativo} are satisfied. The latter result provides the following properties.
    \smallskip\textit{1. Spectral decomposition}: there exists a neighborhood $I_1 \subset I_0$ of zero such that $\tilde{\cl}_{h,\eps}$ has a simple dominating eigenvalue for each $\eps\in I_1$. More precisely, there exist $\eta_1,\eta_2>0$, $c_1\geq 0$ and $\lambda_\eps\in\C$, $f_\eps\in \Bexp(\R^d)$, $\mu_\eps\in (\Bexp(\R^d))^\ast$ and a bounded linear operator $Q_\eps:\Bexp(\R^d)\to\Bexp(\R^d)$ such that 
    \begin{equation}\label{eq: Decomposizione Spettrale}
        \tilde{\cl}_{h,\eps}=\lambda_\eps(f_\eps\otimes\mu_\eps)+Q_\eps
    \end{equation}
    with the following properties:
    \begin{equation}\label{eq: Autovalori e Automisure}
        \tilde{\cl}_{h,\eps}\, f_\eps=\lambda_\eps f_\eps\ ,\qquad \mu_\eps\, \tilde{\cl}_{h,\eps}=\lambda_\eps\mu_\eps\ ,\qquad \mu_\eps\, Q_\eps=0\ ,\qquad Q_\eps f_\eps=0\ ,\qquad \|Q_\eps^n\|_{\Ll(\Bexp)}\leq c_1(1-\eta_1-\eta_2)^n\ .
    \end{equation}
    Moreover, it is possible to require the following normalization conditions\footnote{The second condition holds by Theorem \ref{teo: TeoPerturbativo} $(iii)$. For the first condition, note that $\mu_0$ corresponds to integration with respect to Lebesgue measure; since $f_\eps\to f_0$ in $\Bexp(\R^d)$ as $\eps\to 0$, and $\langle\mu_0,f_0\rangle=1$, by continuity $\langle\mu_0,f_\eps\rangle\neq 0$ for all $\eps\in I_1$, after shrinking $I_1$, if necessary. Replacing $f_\eps$ by $\frac{f_\eps}{\langle\mu_0,f_\eps\rangle}$ and $\mu_\eps$ by $\langle\mu_0,f_\eps\rangle\, \mu_\eps$, we get $\langle\mu_0,f_\eps\rangle=1$ without changing $\langle\mu_\eps,f_\eps\rangle=1$, the spectral projection $f_\eps\otimes\mu_\eps$, or the eigenvalue $\lambda_\eps$.}: \begin{equation}\label{eq: Normalizzazioni}
        \langle \mu_0,f_\eps\rangle=1\ ,\qquad \langle\mu_\eps, f_\eps\rangle =1\ .
    \end{equation}

    \smallskip \textit{2. Eigenmeasure continuity}: the eigenmeasure map $I_1\ni\eps \mapsto \mu_\eps\in (\Bexp(\R^d))'$ is continuous in the dual norm topology. Therefore, evaluated against the fixed invariant density $f_0\in \Bexp(\R^d)$, we obtain 
    \begin{equation}\label{eq: Continuità Automisura}
        \lim_{\eps \to 0}\, \langle \mu_\eps, f_0 \rangle = \langle \mu_0, f_0 \rangle = \int_{\R^d} f_0(x)\, dx = 1\ .
    \end{equation}
    
    \smallskip \textit{3. Eigenvalue expansion}: the eigenvalue map $I_1\ni \eps \mapsto \lambda_\eps\in \C$ is of class $C^1$. By Corollary \ref{cor:derivative}, Proposition \ref{prop: DerivataTransOp Weak} and the Markov property of $S_h$, we deduce that its derivative at zero is given by 
    \[
        \lambda'(0):=\left(\frac{\de }{\de\eps}\lambda{_\eps}\right)_{|_{\eps=0}}=\langle\mu_0, \tilde{\cl}'_h(0)\, f_0\rangle= -\tau\ .
    \]
    Consequently, the eigenvalue admits the first-order expansion 
    \begin{equation}\label{eq: Sviluppo Autovalore}
        \lambda_\eps = 1 - \tau\eps + o(\eps)\ ,
    \end{equation}
    as $\eps\to 0$.
\end{oss}

We are finally able to prove the main theorem. We state it again for the reader's convenience.

\mymaintheorem*

\begin{myproof}
    To establish this limit, we adapt the proof of \cite[Theorem 2]{EXT} to our current functional setting. Using Corollary \ref{cor:EtoTwistedP} with $\phi\equiv 1$ together with the duality relation \eqref{eq:DualitaOperatoriPerturbati},  we rewrite the probability of not exceeding the threshold $u_n$ up to time $t_{n-1}$ in an operator form as:
     \[
        \mathbb{P}\otimes\mu\left(\max_{k=0,\dots,n-1}g_{x_0}(X_{t_k}^x) \leq u_n\right) = \int_{\R^d}\mathbb{E}\left[\mathds{1}_{B_n^c} (X_{t_{0}}^x)\cdots \mathds{1}_{B_n^c}(X_{t_{n-1}}^x)\right] d\mu(x)=\int_{\R^d}(\cl_{h,n}^n\, f_0)(x)\, dx\ .
    \]
    As noted in Remark \ref{oss: Passaggio da famiglia continua a discreta}, we  recover the discrete family $(\cl_{h,n})_{n\geq 1}$ via $\cl_{h,n} = \tilde{\cl}_{h,\eps_n}$, where $\eps_n$ is defined in \eqref{eq: Eps_n}. Since $\eps_n\to 0$ as $n\to+\infty$, we have $\eps_n\in I_1$ for all sufficiently large $n$, and applying the spectral decomposition \eqref{eq: Decomposizione Spettrale}, we have:
    \[
        \tilde{\cl}_{h,\eps_n}^n\, f_0 = \left(\lambda_{\eps_n}^n(f_{\eps_n}\otimes\mu_{\eps_n})+Q_{\eps_n}^n\right)f_0= \lambda_{\eps_n}^n \langle \mu_{\eps_n}, f_0 \rangle\, f_{\eps_n} + Q_{\eps_n}^n f_0 \ .
    \]
    Integrating over $\R^d$ and applying the normalization $\langle \mu_0, f_{\eps_n} \rangle = 1$, we obtain:
    \[
        \int_{\R^d} (\tilde{\cl}_{h,\eps_n}^n\, f_0)(x) \, dx = \lambda_{\eps_n}^n  \langle\mu_{\eps_n},f_0\rangle + \int_{\R^d}(Q_{\eps_n}^n f_0)(x)\, dx\ .
    \]
    By \eqref{eq: Sviluppo Autovalore}, the leading eigenvalue expands as $\lambda_{\eps_n} = 1 - \tau \eps_n + o(\eps_n)$ as $n\to+\infty$. Because $f_0$ is continuous at $x_0$, the Lebesgue differentiation theorem guarantees that the measure of the shrinking balls scales as 
    \[
        \mu(B_n) = f_0(x_0)\, m(B_n) + o(m(B_n)) = \tau \eps_n + o(\eps_n)\ .
    \]
    Thus we can rewrite the expansion as $n\to+\infty$ in terms of the invariant measure:
    \[
        \lambda_{\eps_n} = 1 - \mu(B_n) + o(\mu(B_n))\ .
    \]
    Since the threshold sequence $(u_n)_{n\in\N}\subset\R^+$ satisfies $n\, \mu(B_n) \to \tau$, we have:
    \[
        \lambda_{\eps_n}^n = \left(1 - \frac{\tau + o(1)}{n}\right)^n \to e^{-\tau}\ ,
    \]
    as $n\to+\infty$.
    For the remaining terms in the spectral decomposition, \eqref{eq: Continuità Automisura} establishes that $\langle\mu_{\eps_n},f_0\rangle \to 1$ as $n\to+\infty$. Furthermore, the uniform spectral gap bound $Q_{\eps_n}^n$ in \eqref{eq: Autovalori e Automisure} forces the operator remainder to vanish:
    \[
        \left|\int_{\R^d}(Q^n_{\eps_n} f_0)(x)\, dx\right|  \leq \|Q_{\eps_n}^n\|_{\Ll(\Bexp)} \|f_0\|_{\Bexp}\leq c_1(1-\eta_1-\eta_2)^n\, \|f_0\|_{\Bexp} \to 0\ ,
    \]
    as $n\to+\infty$.
    Putting all these limits together yields the final statement and the proof is done.
\end{myproof}

\section{Conclusions}\label{sec: conclusion}

We have established an extreme value law for discrete-time
observations of stochastic differential equations with bounded, globally Lipschitz
drifts satisfying an inward-pointing condition outside
a bounded region.
For stationary initial data, the limiting law has extremal
index one for every fixed sampling interval $h>0$.

The main methodological contribution is a regular functional
framework in which the smoothing action of the noise allows shrinking-target perturbations to be treated by classical spectral perturbation theory.
On a Banach space combining $C_b^1$ regularity with
exponentially weighted integrability, we establish a spectral
gap for the unperturbed transfer operator and continuous
differentiability in operator norm of the family of killed
operators parameterized by target volume.
The Kantorovich--Rubinstein norm provides an intermediate
topology for describing the derivatives of the removed mass,
which the diffusion maps into the chosen space of regular
densities.
The resulting first-order expansion of the leading eigenvalue,
together with the spectral decomposition, yields the extreme
value law without a separate verification of probabilistic
mixing conditions.

%\subsection*{Further developments}
%\phantomsection
%\addcontentsline{toc}{subsection}{Further developments}

A natural direction is to extend this construction to
diffusions with multiplicative noise and to processes driven
by suitable L\'evy noise.
The central questions are whether the transition operators
provide sufficient regularization and whether suitable
Lyapunov estimates yield the required control at infinity.
These properties would guide the choice of functional spaces
and the analysis of shrinking-target perturbations.
Processes driven by fractional Brownian motion present
an additional difficulty: the memory of the noise generally
requires a formulation beyond a Markov semigroup on the
original state space.

Another direction concerns hitting-time statistics for
non-autonomous SDEs, motivated in particular by applications
involving time-dependent models
\cite{faranda2024statistical,galatolo2026mathematical}.
When the drift depends on time, both the evolving distribution
and the probabilities of occupying a target may vary along
the observation sequence.
A central problem is to identify an appropriate normalization,
potentially based on cumulative probabilities of occupying
the targets, and conditions ensuring an exponential or another
non-degenerate limiting law.
From the operator perspective, this requires replacing the
iteration of a single transfer operator by compositions
of time-dependent operators and controlling the corresponding
loss of memory and rare-event perturbations.

\begin{appendix}

\section{Mixed derivative estimates of the heat kernel with constant diffusion}\label{app: StimaD derivata mista}

In this appendix we establish a mixed-derivative estimate for multidimensional stochastic differential equations driven by Brownian motion. We focus specifically on the regime where the drift coefficient $b(t,x)$ is measurable in $t$ and globally Lipschitz in $x$ (uniformly in $t$), in particular, it may be unbounded with at most linear growth and the diffusion matrix is the identity. To overcome the lack of smoothness in the drift, we adopt the strategy developed in \cite{MPZ}, which introduces a deterministic regularized flow along which the spatial coefficients are frozen. This technique ensures that the resulting constants depend solely on the structural parameters, including the Lipschitz bound, and not on higher derivatives of the mollified drift. Building upon this framework, we derive a mixed-derivative bound that was not present in the original article.

\smallskip In the following, we denote by $\langle\cdot,\cdot\rangle$ and $|\cdot|$ the Euclidean scalar product and norm on $\R^d$. We use the notation $\nabla$ and $\nabla^2$ to denote respectively the gradient and Hessian matrix of real valued smooth functions on $\R^d$. By extension, we denote by $\nabla^j$ the $j$-th order derivative of a smooth function. We use the letter $C$ for constants that may depend on the parameters $\Theta=(T,\beta,\kappa_1,d)$. We reserve the notation $c$ for constants which only depend on the quantities $\beta,\kappa_1$ and $d$ but not on $T$. Both types of constants might change from line to line. Other possible dependencies are explicitly indicated when needed. Finally, we frequently use the notation $\lesssim$: for two quantities $Q_1$ and $Q_2$, we mean by $Q_1\lesssim Q_2$ that there exists $C$ with the same previously described dependence such that $Q_1\leq C\, Q_2$.

\smallskip Let us consider the non-degenerate diffusion process\footnote{In \cite{MPZ}, the non-degenerate diffusion process is $dX_t=b(t,X_t)\, dt+\sigma(t,X_t)\, dW_t$, where the diffusion coefficient satisfies hypothesis $(H_{\alpha}^\sigma)$ (see equations (1.5) and (1.6) in the aforementioned article). Here we work with a less general SDE.}
\begin{equation}\label{eq: MPZ 2.3}
    \begin{cases}
        dX_t=b(t,X_t)\, dt+ dW_t\ ,\qquad t\geq 0\\
        X_0=x
    \end{cases}
\end{equation}
where $(W_t)_{t\geq 0}$ is the $d$-dimensional Brownian motion on the classical Wiener space $(\Omega,\mathcal{F},(\mathcal{F}_t)_{t\geq 0},\mathbb{P})$. We assume that the drift $b$ satisfies hypothesis $(H_\beta^b)$, namely,
\begin{description}
\item[$(H_\beta^b)$] there exist a positive constant $\kappa_1>0$ and a parameter $\beta\in(0,1]$ such that for every $x,y\in\R^d$ and $t\geq 0$
    \begin{equation}\label{eq: MPZ 1.7}
        |b(t,0)|\leq\kappa_1\ ,\qquad |b(t,x)-b(t,y)|\leq\kappa_1\max\{|x-y|^\beta,|x-y|\}\ .
    \end{equation}
\end{description}

Let $\Theta=\Theta(T,\beta,\kappa_1,d)$ be the parameter set appearing in \cite[equation (1.14)]{MPZ} and $\theta_{t,s}^{(1)}(x)$ be the deterministic flow solving \cite[equation (1.11)]{MPZ}. For any $T>0$ and $\eps\in[0,T)$, we write
\[
    \mathbb{D}_\eps^T:=\{(s,t)\in [0,+\infty)^2 : \eps<t-s<T\}\ .
\]

\begin{prop}\label{prop: Stima Mixed Derivative}
    Under $(H_1^b)$, for any $T>0$, $(s,t)\in\mathbb{D}_0^T$ and $x\in\R^d$, the unique strong solution $(X_{t,s}^x)_{t\geq s}$ starting from $x$ at time $s$ of \eqref{eq: MPZ 2.3} admits a density $p(s,x,t,y)$ which is continuous in $(x,y)\in\R^d\times\R^d$ and such that, for each fixed $(s,t)\in\mathbb{D}_0^T$ and $y\in\R^d$, the mixed derivative $\nabla_x\nabla_y p(s,x,t,y)$ exists and is continuous in $x$. Moreover, there exist  constants $\lambda_3,C_3>0$ depending on $\Theta$ such that for all $(s,t)\in\mathbb{D}_0^T$ and $x,y\in\R^d$, 
    \[
        |\nabla_x\nabla_y\, p(s,x,t,y)|\leq C_3(t-s)^{-1}\, g_{\lambda_3}(t-s,\theta_{t,s}^{(1)}(x)-y)\ .
    \]

    \begin{myproof}
        The existence and continuity of the density $p(s,x,t,y)$ follows from \cite[Theorem 1.2]{MPZ}. We therefore focus on the mixed-derivative estimate.
        We write
        \[
            \mathbf{g}_{\lambda}(s,x,t,y):=g_\lambda(t-s,\theta_{t,s}^{(1)}(x)-y)=(t-s)^{-\frac{d}{2}}\, e^{-\frac{\lambda}{t-s}|\theta_{t,s}^{(1)}(x)-y|^2}
        \]
        Let us fix $\eps\in (0,1]$. $p_\eps(s,x,t,y)$ be the heat kernel of the solution $X_{t,s}^{\eps}(x)$ of the system with smooth coefficient corresponding to \eqref{eq: MPZ 2.3}, i.e., with the mollified coefficient $b_\eps$ in the place of $b$ ($\sigma=I$ is already smooth, so it does not need to be regularized). Fixing the \textit{freezing parameters} $(\tau,\xi)\in\R^+\times\R^d$, let us denote by $\tilde{p}_\eps^{(\tau,\xi)}(s,x,t,y)$ the Gaussian density of the process $\tilde{X}_{t,s}^{\eps^{ (\tau,\xi)}}(x)$ (see \cite[equation (2.8)]{MPZ} and note that $\mathcal{C}_{t,s}^{(\tau,\xi)}=(t-s)I$ because $\sigma=I$). If $\theta_{t,s}^{\eps}$ denotes the deterministic flow associated with the mollified drift $b_\eps$, by \cite[Theorem 1.2 and Lemmas 1.1, 2.4]{MPZ} applied to $b_\eps$, all intermediate Gaussian estimates for $p_\eps$ and $\tilde{p}_\eps^{(\tau,\xi)}$ below hold with the $\eps$-independent flow $\theta_{t,s}^{(1)}$ with constants depending only on $\Theta$.
        
        By the \textit{Duhamel-type representation formula in the context of smooth coefficient} (see \cite[equation (2.12)]{MPZ}), recalling that the operator  $\mathcal{A}_{r,z}^{(\tau,\xi)}=B_{r,z}^{(\tau,\xi)}\cdot\nabla_z=(b_\eps(r,z)-b_\eps(r,\theta^\eps_{r,\tau}(\xi)))\cdot\nabla_z$, we have
        \[
            p_\eps(s,x,t,y)=\tilde{p}_\eps^{(\tau,\xi)}(s,x,t,y)+\int_s^t\int_{\R^d}p_\eps(s,x,r,z)\, \mathcal{A}_{r,z}^{(\tau,\xi)}\, \tilde{p}_\eps^{(\tau,\xi)}(r,z,t,y)\, dz\, dr\ .
        \]
        To justify differentiating under the integral sign in the Duhamel representation for fixed $\eps>0$ and arbitrary freezing parameters $(\tau,\xi)\in\R^+\times\R^d$, recall that $b_\eps$ is smooth with bounded spatial derivatives of all orders, thus $p_\eps(s,x,r,z)$ and its spatial derivatives satisfy Gaussian bounds (with constants allowed to depend on $\eps$). First, differentiating in $x$ is justified by dominated convergence, because $|\nabla_x\, p_\eps(s,x,r,z)|$ has an integrable $(r-s)^{-\frac{1}{2}}$ singularity at $r=s$ and $|\nabla_z\, \tilde{p}_\eps^{(\tau,\xi)}(r,z,t,y)|$ has an integrable $(t-r)^{-\frac{1}{2}}$ singularity at $r=t$. Next, to differentiate the resulting integral with respect to $y_j$ (for each component $\de_{x_i}p_\eps(s,x,r,z)$), we split the time integral at $m:=\frac{s+t}{2}$. On $[s,m]$, the factor $\nabla_z\, \de_{y_j}\tilde{p}_\eps^{(\tau,\xi)}(r,z,t,y)$ is smooth and uniformly bounded with Gaussian tails in $z$ because $t-r\geq \frac{t-s}{2}>0$. On $[m,t]$, we have $r-s\geq\frac{t-s}{2}>0$, so for fixed $\eps>0$ the map $z\mapsto\de_{x_i}p_\eps(s,x,r,z)B_{r,z}^{(\tau,\xi)}$ is smooth with uniformly bounded $z$-derivative and Gaussian decay in $z$; integrating by parts in $z$ gives 
        \[
            -\int_{\R^d}\mathrm{div}_z\left(\de_{x_i} p_\eps(s,x,r,z)\, B_{r,z}^{(\tau,\xi)}\right) \tilde{p}_{\eps}^{(\tau,\xi)}(r,z,t,y)\, dz\ .
        \]
        Differentiating in $y_j$ under the integration in $dr$ over $[m,t]$ is now justified by dominated convergence because $\de_{y_j}\tilde{p}_\eps^{(\tau,\xi)}(r,z,t,y)$ has only an integrable $(t-r)^{-\frac{1}{2}}$ singularity at $r=t$. Integrating by parts in $z$ once more for each fixed $r\in (m,t)$ gives
        \begin{align*}
             \nabla_x\nabla_y\, p_\eps(s,x,t,y)
             &=\nabla_x\nabla_y\, \tilde{p}_\eps^{(\tau,\xi)}(s,x,t,y)+\int_s^t\int_{\R^d}\nabla_x\, p_\eps(s,x,r,z)\, B_{r,z}^{(\tau,\xi)}\cdot \nabla_z\nabla_y\, \tilde{p}_\eps^{(\tau,\xi)}(r,z,t,y)\, dz\, dr\\
            &=-\nabla_x^2\, \tilde{p}_\eps^{(\tau,\xi)}(s,x,t,y)-\int_s^t\int_{\R^d}\nabla_x\, p_\eps(s,x,r,z)\, B_{r,z}^{(\tau,\xi)}\cdot \nabla_z^2\, \tilde{p}_\eps^{(\tau,\xi)}(r,z,t,y)\, dz\, dr\ ,
        \end{align*}
         where in the last equality we use $\nabla_x\, \tilde{p}_\eps^{(\tau,\xi)}(s,x,t,y)=-\nabla_y\, \tilde{p}_\eps^{(\tau,\xi)}(s,x,t,y)$ and $\nabla_y\, \tilde{p}_\eps^{(\tau,\xi)}(r,z,t,y)=-\nabla_z\, \tilde{p}_\eps^{(\tau,\xi)}(r,z,t,y)$ for the first addend and for the second factor inside the integral. Taking $(\tau,\xi)=(t,y)$,  we get
        \[
            \nabla_x\nabla_y\, p_\eps(s,x,t,y)=-\nabla_x^2\, \tilde{p}_\eps^{(t,y)}(s,x,t,y)-\int_s^t\int_{\R^d}\nabla_x\, p_\eps(s,x,r,z)\, B_{r,z}^{(t,y)}\cdot \nabla_z^2\, \tilde{p}_\eps^{(t,y)}(r,z,t,y)\, dz\, dr\ .
        \]
        Hence
        \[
            |\nabla_x\nabla_y\, p_\eps(s,x,t,y)|\lesssim (t-s)^{-1}\, \mathbf{g}_{\tilde{\lambda}_2}(s,x,t,y)+\int_s^t\int_{\R^d}(r-s)^{-\frac{1}{2}}\, \mathbf{g}_{\lambda_1}(s,x,r,z)|z-\theta_{r,t}^\eps(y)|(t-r)^{-1}\, \mathbf{g}_{\tilde{\lambda}_2} (r,z,t,y)\, dz\, dr\ ,
        \]
        where we applied \cite[Lemma 2.4]{MPZ} twice with $\alpha=0$, \cite[Theorem 1.2]{MPZ} and the Lipschitz continuity of $b_\eps$. By \cite[Lemma 1.1]{MPZ}, we get
        \begin{equation}\label{eq: Singolarita Ridotta Riuso}
            \frac{|z-\theta_{r,t}^\eps(y)|}{(t-r)^{\frac{1}{2}}}\, \mathbf{g}_{\tilde{\lambda}_2}(r,z,t,y)\lesssim\frac{|y-\theta_{t,r}^{(1)}(z)|+(t-r)^{\frac{1}{2}}}{(t-r)^{\frac{1}{2}}}\, \mathbf{g}_{\tilde{\lambda}_2}(r,z,t,y)\lesssim \mathbf{g}_{\tilde{\lambda}'_2}(r,z,t,y)\ .
        \end{equation}
        Therefore
        \begin{align*}
            |\nabla_x\nabla_y\, p_\eps(s,x,t,y)|&\lesssim  (t-s)^{-1}\, \mathbf{g}_{\tilde{\lambda}_2}(s,x,t,y)+\int_s^t\int_{\R^d}(r-s)^{-\frac{1}{2}}\, \mathbf{g}_{\lambda_1}(s,x,r,z)\, (t-r)^{-\frac{1}{2}}\frac{|z-\theta_{r,t}^\eps(y)|}{(t-r)^{\frac{1}{2}}}\, \mathbf{g}_{\tilde{\lambda}_2} (r,z,t,y)\, dz\, dr\ ,\\
            &\lesssim (t-s)^{-1}\, \mathbf{g}_{\tilde{\lambda}_2}(s,x,t,y)+\int_s^t(r-s)^{-\frac{1}{2}}(t-r)^{-\frac{1}{2}}\left(\int_{\R^d}\mathbf{g}_{\lambda_1}(s,x,r,z)\, \mathbf{g}_{\tilde{\lambda}'_2}(r,z,t,y)\, dz\right)dr\\
            &\lesssim (t-s)^{-1}\, \mathbf{g}_{\tilde{\lambda}_2}(s,x,t,y)+\left(\int_s^t(r-s)^{-\frac{1}{2}}(t-r)^{-\frac{1}{2}}\, dr\right) \mathbf{g}_{\lambda_3}(s,x,t,y)\\
            &\lesssim (t-s)^{-1}\, \mathbf{g}_{\lambda_3}(s,x,t,y)\ ,
        \end{align*}
        where in the third inequality we applied  \cite[Lemma 2.5]{MPZ} to bound the Gaussian convolution and in the fourth one we use the fact that 
        \[
            \int_s^t(r-s)^{-\frac{1}{2}}(t-r)^{-\frac{1}{2}}\, dr=B\left(\frac{1}{2},\frac{1}{2}\right)\leq B\left(\frac{1}{2},\frac{1}{2}\right)\, \frac{T}{t-s}\ ,
        \]
        where we recall that $B$ denotes the $\beta$-function and the equality is obtained via the substitution $u=\frac{r-s}{t-s}$, while the inequality holds
        because $T\geq t-s$, thus even in the second addend we recover the singularity factor $(t-s)^{-1}$. Hence
        \begin{equation}\label{eq: StimaMisteKernelIntermezzo}
             |\nabla_x\nabla_y\, p_\eps(s,x,t,y)|\leq C_3(t-s)^{-1}\, g_{\lambda_3}(t-s,\theta_{t,s}^{(1)}(x)-y)\ .
        \end{equation}

       By inspecting the constants involved in \cite[Lemmas 1.1, 2.4 and Theorem 1.2]{MPZ}, we note that $\lambda_3$ and $C_3$ are actually \textit{independent} of the continuity and of the derivatives of the mollified drift $b_\eps$ or of the mollification parameter $\eps$, but only depend on the parameters gathered in $\Theta$. 

       \smallskip To conclude, it is sufficient to verify the sensitivity with respect to the variable $x$ (for each fixed $s<t$ and $y\in\R^d$): indeed, since $\nabla_y\, p_\eps(s,x,t,y)\to\nabla_y\, p(s,x,t,y)$ locally uniformly by \cite[Lemma A.1]{MPZ} together with \cite[Remark 1.3]{MPZ}, uniform boundedness \eqref{eq: StimaMisteKernelIntermezzo} and equicontinuity of $x\mapsto\nabla_x\nabla_y\, p_\eps(s,x,t,y)$ for fixed $y$ imply via Ascoli--Arzelà theorem that $\nabla_x\nabla_y\, p(s,x,t,y)$ exists, is continuous in $x$ and satisfies the bound in the statement. 
       
       Namely, we show the following. Let $\gamma\in (0,1)$. Then there exist constants $\lambda,C>0$ depending only on $\Theta$ such that, for all $(s,t)\in\mathbb{D}_0^T$ and $x,x',y\in\R^d$,
       \[
            |\nabla_x\nabla_y\, p_\eps(s,x,t,y)-\nabla_x\nabla_y\, p_\eps(s,x',t,y)|\lesssim |x-x'|^\gamma (t-s)^{-1-\frac{\gamma}{2}}\left(\mathbf{g}_{\lambda}(s,x,t,y)+\mathbf{g}_\lambda(s,x',t,y)\right)\ .
       \]

         First of all we consider the  ``\textit{off-diagonal case}'': $|x-x'|^2>\frac{t-s}{4}$. By \eqref{eq: StimaMisteKernelIntermezzo}, we clearly have
       \begin{align*}
            |\nabla_x\nabla_y\, p_\eps(s,x,t,y)-\nabla_x\nabla_y\, p_\eps(s,x',t,y)|&\lesssim (t-s)^{-1}\left(\mathbf{g}_{\lambda}(s,x,t,y)+\mathbf{g}_{\lambda}(s,x',t,y)\right)\\
            &\lesssim|x-x'|^\gamma (t-s)^{-1-\frac{\gamma}{2}}\left(\mathbf{g}_{\lambda}(s,x,t,y)+\mathbf{g}_\lambda(s,x',t,y)\right)\ ,
       \end{align*}
       where the factor multiplying the Gaussians is due to the fact that in this case
       \[
            1>\frac{(t-s)^{\frac{1}{2}}}{2\, |x-x'|}\Longleftrightarrow 1<\frac{2\, |x-x'|}{(t-s)^{\frac{1}{2}}}\Longleftrightarrow 1\leq \left(\frac{2\, |x-x'|}{(t-s)^{\frac{1}{2}}}\right)^\gamma
       \]
       for each $\gamma\in(0,1)$.
       
       Next we restrict to the so-called ``\textit{diagonal case}'': $|x-x'|^2\leq\frac{t-s}{4}$.
      By the Duhamel representation formula (see \cite[equation (2.12)]{MPZ}), we have 
      \begin{align}\label{eq: DaDerivareDueVolte}
            \begin{split}p_\eps(s,x,t,y)-p_\eps(s,x',t,y)=\ &\tilde{p}_\eps^{(\tau,\xi)}(s,x,t,y)-\tilde{p}_\eps^{(\tau,\xi)}(s,x',t,y)\\
            &+\int_s^t\int_{\R^d}(p_\eps(s,x,r,z)-p_\eps(s,x',r,z))\, \mathcal{A}_{r,z}^{(\tau,\xi)}\, \tilde{p}_\eps^{(\tau,\xi)}(r,z,t,y)\, dz\, dr\ .
            \end{split}
        \end{align}
        Hence, differentiating first with respect to $x$ and $y$, by arbitrary freezing parameters $(\tau,\xi)\in \R^+\times\R^d$, as done in the previous differentiation under integral sign, and then taking $\tau=t$ and $\xi=y$, we obtain
        \[
            \nabla_x\nabla_y\, p_\eps(s,x,t,y)-\nabla_x\nabla_y\, p_\eps(s,x',t,y)=I_1(s,t,x,x',y)+I_2(s,t,x,x',y)
        \] 
        for some $I_i(s,t,x,x',y)$ with $i=1,2$ that we explicitly compute and then estimate below. 
        
        Observe that, for any freezing couple $(\tau,\xi)\in \R^+\times\R^d$ and $h\in\R^d$,
        \begin{equation}\label{eq: PassaggioDerivateSecondeXY}
                \qquad \nabla_x^2\, \tilde{p}_\eps^{(\tau,\xi)}(s,x+h,t,y)=\nabla_y^2\, \tilde{p}_\eps^{(\tau,\xi)}(s,x,t,y-h)\ .
        \end{equation}
         Thus differentiating \eqref{eq: DaDerivareDueVolte}, using that $\nabla_x\, \tilde{p}_\eps^{(\tau,\xi)}(s,x,t,y)=-\nabla_y\, \tilde{p}_\eps^{(\tau,\xi)}(s,x,t,y)$ and \eqref{eq: PassaggioDerivateSecondeXY} applied with $h=x-x'$,  we have
        \begin{align*}
            \nabla_x\nabla_y\, \tilde{p}_\eps^{(\tau,\xi)}(s,x,t,y)-\nabla_x\nabla_y\, \tilde{p}_\eps^{(\tau,\xi)}(s,x',t,y)&=-\nabla_x^2\, \tilde{p}_\eps^{(\tau,\xi)}(s,x,t,y)+\nabla_y^2\, \tilde{p}_\eps^{(\tau,\xi)}(s,x',t,y)\\
            &=-\nabla_y^2\, \tilde{p}_\eps^{(\tau,\xi)}(s,x',t,y+x'-x)+\nabla_y^2\, \tilde{p}_\eps^{(\tau,\xi)}(s,x',t,y)\ .
        \end{align*}
        Evaluating in $\tau=t$ and $\xi=y$, we get
        \[
            I_1(s,t,x,x',y)=-\nabla_y^2\, \tilde{p}_\eps^{(t,y)}(s,x',t,y+x'-x)+\nabla_y^2\, \tilde{p}_\eps^{(t,y)}(s,x',t,y)\ ,
        \]
        where $\nabla_y^j\, \tilde{p}_\eps^{(t,y)}(s,x,t,w)$ denotes the $j$-th derivative with respect to the fourth argument $w$, keeping the freezing parameters $(t,y)$ fixed.

        Similarly to what is done in the proof of [Lemma A.1] (see \cite[equation (A.4)]{MPZ}), by \cite[Lemma 2.3]{MPZ}, for $j\in\N$ and $h\in\R^d$ with $|h|^2\leq \frac{t-s}{4}$,
         \begin{align}\label{eq: MPZ A.4}
           \begin{split}\left|\nabla_y^j\, \tilde{p}_\eps^{(t,y)}(s,x,t,y+h)-\nabla_y^j\, \tilde{p}_\eps^{(t,y)}(s,x,t,y)\right|&\leq |h|\sup_{\rho\in [0,1]}\left|\nabla_y^{j+1}\, \tilde{p}_\eps^{(t,y)}(s,x,t,y+\rho h)\right|\\
           &\lesssim |h|(t-s)^{-\frac{j+1}{2}}\sup_{\rho\in[0,1]}\, g_{\lambda_1}(t-s,\theta_{t,s}^{(1)}(x)-(y+\rho h))\\
           &\lesssim |h|(t-s)^{-\frac{j+1}{2}}\, g_{\lambda_1}(t-s,\theta_{t,s}^{(1)}(x)-y)\ ,
           \end{split}
         \end{align}
         using \cite[Lemma 1.1]{MPZ} for the last step, after shrinking $\lambda_1$ to a smaller parameter, still denoted by $\lambda_1$ by a small abuse of notation.
         On the other hand, we also have
         \[
            \left|\nabla_y^j\, \tilde{p}_\eps^{(t,y)}(s,x,t,y+h)\right|\lesssim (t-s)^{-\frac{j}{2}}\, g_{\lambda_1}(t-s,\theta_{t,s}^{(1)}(x)-y)\ .
         \]
         Thus, by interpolation, we get, for any $\gamma\in (0,1)$,
        \[
             \left|\nabla_y^j\, \tilde{p}_\eps^{(t,y)}(s,x,t,y+h)-\nabla_y^j\, \tilde{p}_\eps^{(t,y)}(s,x,t,y)\right|\lesssim |h|^\gamma(t-s)^{-\frac{j+\gamma}{2}}\, g_{\lambda_1}(t-s,\theta_{t,s}^{(1)}(x)-y)\ .
        \]
        Hence, we can bound the first term as follows:
         \[
            |I_1(s,t,x,x',y)|\lesssim |x-x'|^\gamma (t-s)^{-1-\frac{\gamma}{2}}\, \mathbf{g}_{\lambda_1}(s,x',t,y)\leq |x-x'|^\gamma (t-s)^{-1-\frac{\gamma}{2}}\left(\mathbf{g}_{\lambda_1}(s,x,t,y) +\mathbf{g}_{\lambda_1}(s,x',t,y)\right)\ .
         \]
        
        To treat the remaining term, recalling that $\mathcal{A}_{r,z}^{(\tau,\xi)}=B_{r,z}^{(\tau,\xi)}\cdot\nabla_z=(b_\eps(r,z)-b_\eps(r,\theta_{r,\tau}^{\eps}(\xi)))\cdot\nabla_z$, differentiating  \eqref{eq: DaDerivareDueVolte} (by a small abuse of notation, we still denote by $\nabla_x$ the derivative with respect to the second argument of the kernel $p$ evaluated at those points), we get 
        \[
            \int_s^t\int_{\R^d}\nabla_x\left(p_\eps(s,x,r,z)-p_\eps(s,x',r,z)\right) B_{r,z}^{(\tau,\xi)}\cdot\nabla_z\nabla_y\, \tilde{p}_\eps^{(\tau,\xi)}(r,z,t,y)\, dz\, dr\ .
        \]
        Using that $\nabla_y\, \tilde{p}_\eps^{(\tau,\xi)}(r,z,t,y)=-\nabla_z\, \tilde{p}_\eps^{(\tau,\xi)}(r,z,t,y)$ for the second factor inside the integral, we can rewrite the above integral as follows
        \[
             \int_s^t\int_{\R^d}\nabla_x\left(p_\eps(s,x,r,z)-p_\eps(s,x',r,z)\right) B_{r,z}^{(\tau,\xi)}\cdot\left(-\nabla_z^2\, \tilde{p}_\eps^{(\tau,\xi)}(r,z,t,y)\right) dz\, dr\ .
        \]
        Evaluating in $\tau=t$ and $\xi=y$, we get
        \[
            I_2(s,t,x,x',y)=-\int_s^t\int_{\R^d}\nabla_x\left(p_\eps(s,x,r,z)-p_\eps(s,x',r,z)\right) B_{r,z}^{(t,y)}\cdot\nabla_z^2\, \tilde{p}_\eps^{(t,y)}(r,z,t,y)\, dz\, dr\ .
        \]
        It holds (see \cite[Lemma A.1 point $(C_1)$]{MPZ}) 
        \[
            \left|\nabla_x\, p_\eps(s,x,r,z)-\nabla_x\, p_\eps(s,x',r,z)\right|\lesssim |x-x'|^\gamma (r-s)^{-\frac{1+\gamma}{2}}\left(\mathbf{g}_{\lambda_2}(s,x,r,z)+\mathbf{g}_{\lambda_2}(s,x',r,z)\right)
        \]
         for $\gamma\in (0,1)$. 
         Therefore the previous bound combined with \cite[Lemma 2.4]{MPZ}, applied with $\alpha=0$, gives that
         \begin{align*}
             |I_2(s,t,x,x',y)|&\lesssim |x-x'|^\gamma\int_s^t\int_{\R^d}(r-s)^{-\frac{1+\gamma}{2}}\left(\mathbf{g}_{\lambda_2}(s,x,r,z)+\mathbf{g}_{\lambda_2}(s,x',r,z)\right)|z-\theta_{r,t}^\eps(y)|\, (t-r)^{-1}\, \mathbf{g}_{\lambda_2}(r,z,t,y)\, dz\, dr\\
             &\leq|x-x'|^\gamma\int_s^t\int_{\R^d}(r-s)^{-\frac{1+\gamma}{2}}\left(\mathbf{g}_{\lambda_2}(s,x,r,z)+\mathbf{g}_{\lambda_2}(s,x',r,z)\right)(t-r)^{-\frac{1}{2}}\frac{|z-\theta_{r,t}^\eps(y)|}{(t-r)^{\frac{1}{2}}}\, \mathbf{g}_{\lambda_2}(r,z,t,y)\, dz\, dr\\
             &\lesssim |x-x'|^\gamma\int_s^t(r-s)^{-\frac{1+\gamma}{2}}(t-r)^{-\frac{1}{2}}\int_{\R^d}\left(\mathbf{g}_{\lambda_2}(s,x,r,z)+\mathbf{g}_{\lambda_2}(s,x',r,z)\right)\, \mathbf{g}_{\lambda_2'}(r,z,t,y)\, dz\, dr\\
             &\lesssim |x-x'|^\gamma (t-s)^{-\frac{\gamma}{2}}\left(\mathbf{g}_{\lambda_2''}(s,x,t,y)+\mathbf{g}_{\lambda_2''}(s,x',t,y)\right)\\
             &\leq T\, |x-x'|^\gamma (t-s)^{-1-\frac{\gamma}{2}}\left(\mathbf{g}_{\lambda_2''}(s,x,t,y)+\mathbf{g}_{\lambda_2''}(s,x',t,y)\right)\ ,
         \end{align*}
         where in the fourth inequality we use the fact that, via the substitution $u=\frac{r-s}{t-s}$, we get
        \[
            \int_s^t(r-s)^{-\frac{1+\gamma}{2}}(t-r)^{-\frac{1}{2}}\, dr=(t-s)^{-\frac{\gamma}{2}}B\left(\frac{1-\gamma}{2},\frac{1}{2}\right)\ .
        \]
        
         Hence, combining the above calculations, we get the existence of the desired constants $C=C(\Theta)$ and $\lambda=\lambda(\Theta)$ such that 
         \[
            |I_1(s,t,x,x',y)+I_2(s,t,x,x',y)|\lesssim_C |x-x'|^\gamma(t-s)^{-\left(1+\frac{\gamma}{2}\right)}\left(\mathbf{g}_{\lambda}(s,x,t,y)+\mathbf{g}_{\lambda}(s,x',t,y)\right)\ ,
         \]
         which concludes the proof of the equicontinuity of $\nabla_x\nabla_y\, p_\eps(s,x,t,y)$ and therefore the proof of the proposition.
    \end{myproof}
\end{prop}

\end{appendix}

\section*{Acknowledgements} 
We acknowledge the MIUR Excellence Department Project awarded to the
Department of Mathematics, University of Pisa, CUP I57G22000700001.
 The authors would like to express their gratitude to Professor Antonello Pesce for his valuable assistance and suggestions regarding the proof of the mixed-derivative estimate presented in Appendix \ref{app: StimaD derivata mista}.

\section*{AI disclosure}

AI was used primarily for proofreading the text and for improving its presentation and style. It proved helpful in optimizing and verifying some computations. All the ideas contained in this work are to be attributed to the authors.

\section*{Competing interests}

    The authors have no conflict of interests or competing interests to declare that are relevant to the content of this article.

\bibliographystyle{alpha}\bibliography{bib}

@article{EXT,
  title={Extreme value theory and {P}oisson statistics for discrete time samplings of stochastic differential equations},
  author={Flandoli, Franco and Galatolo, Stefano and Giulietti, Paolo and Vaienti, Sandro},
  journal={Communications in Mathematical Physics},
  volume={406},
  number={10},
  pages={231},
  year={2025},
  publisher={Springer}
}

@book{HennionHerve,
  title={Limit theorems for {M}arkov chains and stochastic properties of dynamical systems by quasi-compactness},
  author={Hennion, Hubert and Herv{\'e}, Lo{\"\i}c},
  volume={1766},
  series={Lecture Notes in Mathematics},
  year={2002},
  publisher={Springer}
}

@article{MPZ,
  title={Density and gradient estimates for non degenerate {B}rownian {SDE}s with unbounded measurable drift},
  author={Menozzi, St{\'e}phane and Pesce, Antonello and Zhang, Xicheng},
  journal={Journal of Differential Equations},
  volume={272},
  pages={330--369},
  year={2021},
  publisher={Elsevier}
}

@incollection{Friedman64,
  title={Stochastic differential equations and applications},
  author={Friedman, Avner},
  booktitle={Stochastic Differential Equations},
  pages={75--148},
  year={1975},
  publisher={Springer}
}

@article{KL09,
  title={Rare events, escape rates and quasistationarity: some exact formulae},
  author={Keller, Gerhard and Liverani, Carlangelo},
  journal={Journal of Statistical Physics},
  volume={135},
  number={3},
  pages={519--534},
  year={2009},
  publisher={Springer}
}

@article{KL12,
  title={Rare events, exponential hitting times and extremal indices via spectral perturbation},
  author={Keller, Gerhard},
  journal={Dynamical Systems},
  volume={27},
  number={1},
  pages={11--27},
  year={2012},
  publisher={Taylor \& Francis}
}

@article{HV04,
  author  = {Haydn, Nicolai and Vaienti, Sandro},
  title   = {The limiting distribution and error terms for return times
             of dynamical systems},
  journal = {Discrete and Continuous Dynamical Systems},
  volume  = {10},
  number  = {3},
  pages   = {589--616},
  year    = {2004},
  doi     = {10.3934/dcds.2004.10.589}
}

@article{AFGTV25,
  author  = {Atnip, Jason and Froyland, Gary and
             Gonz{\'a}lez-Tokman, Cecilia and Vaienti, Sandro},
  title   = {Compound {Poisson} statistics for dynamical systems
             via spectral perturbation},
  journal = {Nonlinearity},
  volume  = {38},
  number  = {8},
  pages   = {085023},
  year    = {2025},
  doi     = {10.1088/1361-6544/adfa62}
}

@article{AHV25,
  author  = {Amorim, Lucas and Haydn, Nicolai and Vaienti, Sandro},
  title   = {Compound {Poisson} distributions for random dynamical
             systems using probabilistic approximations},
  journal = {Stochastic Processes and their Applications},
  volume  = {179},
  pages   = {104511},
  year    = {2025},
  doi     = {10.1016/j.spa.2024.104511}
}

@book{Kato,
  title={Perturbation theory for linear operators},
  author={Kato, Tosio},
  volume={132},
  year={1966},
  publisher={Springer}
}

@article{marinescu1950theorie,
  title={Th{\'e}orie ergodique pour des classes d'op{\'e}rations non compl{\`e}tement continues},
  author={Ionescu Tulcea, C. T. and Marinescu, G.},
  journal={Annals of Mathematics},
  volume={52},
  number={1},
  pages={140--147},
  year={1950},
  publisher={JSTOR}
}

@book{MPZRef26,
  title={Multidimensional diffusion processes},
  author={Stroock, Daniel W. and Varadhan, S. R. Srinivasa},
  year={2007},
  publisher={Springer}
}

@book{lucarini2016extremes,
  title={Extremes and recurrence in dynamical systems},
  author={Lucarini, Valerio and Faranda, Davide and Freitas, Ana Cristina Moreira and Freitas, Jorge Milhazes and Holland, Mark and Kuna, Tobias and Nicol, Matthew and Todd, Mike and Vaienti, Sandro},
  year={2016},
  publisher={John Wiley \& Sons}
}

@article{Hitting,
  title={Hitting time statistics and extreme value theory},
  author={Freitas, Ana Cristina Moreira and Freitas, Jorge Milhazes and Todd, Mike},
  journal={Probability Theory and Related Fields},
  volume={147},
  number={3},
  pages={675--710},
  year={2010},
  publisher={Springer}
}

@article{leadbetter1988extremal,
  title={Extremal theory for stochastic processes},
  author={Leadbetter, Malcolm R. and Rootz{\'e}n, Holger},
  journal={The Annals of Probability},
  volume={16},
  number={2},
  pages={431--478},
  year={1988},
  publisher={JSTOR}
}

@book{leadbetterVecchio,
  title={Extremes and related properties of random sequences and processes},
  author={Leadbetter, Malcolm R. and Lindgren, Georg and Rootz{\'e}n, Holger},
  year={2012},
  publisher={Springer Science \& Business Media}
}

@article{hsing1988exceedance,
  title={On the exceedance point process for a stationary sequence},
  author={Hsing, Tailen and H{\"u}sler, J{\"u}rg and Leadbetter, Malcolm Ross},
  journal={Probability Theory and Related Fields},
  volume={78},
  number={1},
  pages={97--112},
  year={1988},
  publisher={Springer}
}

@article{haydn2005hitting,
  title={Hitting and return times in ergodic dynamical systems},
  author={Haydn, Nicolai and Lacroix, Yves and Vaienti, Sandro},
  journal={The Annals of Probability},
  volume={33},
  number={5},
  pages={2043--2050},
  year={2005}
}

@article{abadi2001inequalities,
  title={Inequalities for the occurrence times of rare events in mixing processes. {T}he state of the art},
  author={Abadi, Miguel and Galves, Antonio},
  journal={Markov Processes and Related Fields},
  volume={7},
  number={1},
  pages={97--112},
  year={2001}
}

@article{nicolai2009compound,
  title={The compound {P}oisson distribution and return times in dynamical systems},
  author={Haydn, Nicolai and Vaienti, Sandro},
  journal={Probability Theory and Related Fields},
  volume={144},
  number={3-4},
  pages={517--542},
  year={2009},
  publisher={Springer}
}

@article{hirata1999statistics,
  title={Statistics of return times: a general framework and new applications},
  author={Hirata, Masaki and Saussol, Beno{\^\i}t and Vaienti, Sandro},
  journal={Communications in Mathematical Physics},
  volume={206},
  number={1},
  pages={33--55},
  year={1999},
  publisher={Springer}
}

@article{faranda2024statistical,
  title={Statistical physics and dynamical systems perspectives on geophysical extreme events},
  author={Faranda, Davide and Messori, Gabriele and Alberti, Tommaso and Alvarez-Castro, Carmen and Caby, Th{\'e}ophile and Cavicchia, Leone and Coppola, Erika and Donner, Reik V. and Dubrulle, B{\'e}reng{\`e}re and Galfi, Vera Melinda and others},
  journal={Physical Review E},
  volume={110},
  number={4},
  pages={041001},
  year={2024},
  publisher={APS}
}

@article{galatolo2026mathematical,
  title={A mathematical framework for linear response theory for nonautonomous systems},
  author={Galatolo, Stefano and Lucarini, Valerio},
  journal={arXiv preprint arXiv:2603.19509},
  year={2026}
}

@article{guivarc1988theoremes,
  title={Th{\'e}or{\`e}mes limites pour une classe de cha{\^\i}nes de {M}arkov et applications aux diff{\'e}omorphismes d'{A}nosov},
  author={Guivarc'h, Yves and Hardy, Jean},
  journal={Annales de l'Institut Henri Poincar{\'e}, Probabilit{\'e}s et Statistiques},
  volume={24},
  number={1},
  pages={73--98},
  year={1988}
}

@article{Glynn11,
  author  = {Glynn, Peter W.},
  title   = {On exponential limit laws for hitting times of rare sets for {Harris} chains and processes},
  journal = {Journal of Applied Probability},
  year    = {2011},
  volume  = {48},
  number  = {A},
  pages   = {319--326},
  doi     = {10.1239/jap/1318940474}
}

@article{KL99,
  author  = {Keller, Gerhard and Liverani, Carlangelo},
  title   = {Stability of the spectrum for transfer operators},
  journal = {Annali della Scuola Normale Superiore di Pisa, Classe di Scienze},
  series  = {4},
  year    = {1999},
  volume  = {28},
  number  = {1},
  pages   = {141--152},
  url     = {https://www.numdam.org/item/ASNSP_1999_4_28_1_141_0/}
}

@incollection{HM11,
  author    = {Hairer, Martin and Mattingly, Jonathan C.},
  title     = {Yet another look at {Harris}' ergodic theorem for {Markov} chains},
  booktitle = {Seminar on Stochastic Analysis, Random Fields and Applications VI},
  editor    = {Dalang, Robert C. and Dozzi, Marco and Russo, Francesco},
  series    = {Progress in Probability},
  volume    = {63},
  publisher = {Birkh{\"a}user},
  address   = {Basel},
  year      = {2011},
  pages     = {109--117},
  doi       = {10.1007/978-3-0348-0021-1_7}
}

@article{KM03,
  author  = {Kontoyiannis, Ioannis and Meyn, Sean P.},
  title   = {Spectral theory and limit theorems for geometrically ergodic {Markov} processes},
  journal = {The Annals of Applied Probability},
  year    = {2003},
  volume  = {13},
  number  = {1},
  pages   = {304--362},
  doi     = {10.1214/aoap/1042765670}
}

@article{AbadiSaussol11,
  author  = {Abadi, Miguel and Saussol, Beno{\^i}t},
  title   = {Hitting and returning to rare events for all alpha-mixing processes},
  journal = {Stochastic Processes and their Applications},
  year    = {2011},
  volume  = {121},
  number  = {2},
  pages   = {314--323},
  doi     = {10.1016/j.spa.2010.11.001}
}

@article{HP10,
  author  = {Herv{\'e}, Lo{\"i}c and P{\`e}ne, Fran{\c{c}}oise},
  title   = {The {Nagaev--Guivarc'h} method via the {Keller--Liverani} theorem},
  journal = {Bulletin de la Soci{\'e}t{\'e} Math{\'e}matique de France},
  year    = {2010},
  volume  = {138},
  number  = {3},
  pages   = {415--489},
  doi     = {10.24033/bsmf.2594}
}

@article{Zwe22,
  author  = {Zweim{\"u}ller, Roland},
  title   = {Hitting times and positions in rare events},
  journal = {Annales Henri Lebesgue},
  year    = {2022},
  volume  = {5},
  pages   = {1361--1415},
  doi     = {10.5802/ahl.154}
}
\end{document}